\documentclass[11pt,twoside,reqno]{amsart}
\usepackage[usenames,dvipsnames]{color}
\usepackage{amsmath,amsthm,amsopn,amssymb,a4wide,times,enumerate,varioref}
\usepackage[numbers,sort&compress]{natbib}
\usepackage[mathscr]{eucal} 

\numberwithin{equation}{section}

\newtheorem{theorem}{Theorem}[section]
\newtheorem{lemma}[theorem]{Lemma}
\newtheorem{corollary}[theorem]{Corollary}
\newtheorem{proposition}[theorem]{Proposition}

\newtheorem{definition}[theorem]{Definition}
\newtheoremstyle{note}
     {}
     {}
     {}
     {}
     {\bfseries} 
     {.}
     {.5em}
     {}

\theoremstyle{note}

\newtheorem{remark}[theorem]{Remark}

\newcommand{\ip}[2]{{\left<#1,#2\right>}}
\newcommand{\TT}{\Psi}

\newcommand{\CT}{{\mathfrak{B}}}
\newcommand{\CTs}{{\mathfrak{B}}^*}
\newcommand{\s}{\mathbin{\mbox{{$\>\star\>$}}}}
\newcommand{\st}{\star}
\newcommand{\sh}{\mathbin{\mbox{{$\>\hat{\star}\>$}}}}
\newcommand{\sth}{\hat{\star}}
\newcommand{\sT}{\mathfrak T}
\newcommand{\HP}{H^\infty(\mathcal K_\Psi)}
\newcommand{\HPF}{H^\infty(\mathcal K_\Psi^F)}
\newcommand{\PrF}{\Psi|_{F}}
\newcommand{\leftless}{\le_\ell}
\newcommand{\leftsless}{<_\ell}

\newcommand{\tE}{\tilde{E}}

\begin{document}

\setcounter{page}{1} \title{Interpolation in Semigroupoid Algebras}

\author[M.~A.~Dritschel, S.~A.~M.~Marcantognini and
S.~McCullough]{Michael A. Dritschel$^1$, Stefania Marcantognini$^1$
  and Scott McCullough$^2$}

\address{School of Mathematics and Statistics\\
  Merz Court,\\ University of Newcastle upon Tyne\\
  Newcastle upon Tyne\\
  NE1 7RU\\
  UK}

\email{m.a.dritschel@ncl.ac.uk}

\address{Departamento de Matem\'aticas\\
Instituto Venezolano de Investigaciones Cient\'{\i}ficas\\
Apartado 21827\\
Caracas 1020A\\
Venezuela}

\email{smarcant@ivic.ve}

\address{Department of Mathematics\\
  University of Florida\\
  Box 118105\\
  Gainesville, FL 32611-8105\\
  USA}

\email{sam@math.ufl.edu}

\subjclass[2000]{47A57 (Primary), 47L55, 47L75, 47D25, 47A13, 47B38, 46E22 (Secondary)}

\keywords{interpolation, transfer functions, semigroupoid,
  noncommutative function space, multiply connected domains,
  Nevanlinna-Pick, Carath\'eodory-Fej\'er}

\begin{abstract}
  A seminal result of Agler characterizes the so-called Schur-Agler
  class of functions on the polydisk in terms of a unitary colligation
  transfer function representation.  We generalize this to the unit
  ball of the algebra of multipliers for a family of test functions
  over a broad class of semigroupoids.  There is then an associated
  interpolation theorem.  Besides leading to solutions of the familiar
  Nevanlinna-Pick and Carath\'eodory-Fej\'er interpolation problems
  and their multivariable commutative and noncommutative
  generalizations, this approach also covers more exotic examples.
\end{abstract}

\thanks{${}^1$Research supported by the EPSRC.  \quad ${}^2$Research
  supported by the NSF}


\maketitle


\section{Introduction}

The transfer function realization formalism for contractive
multipliers of (families of) reproducing kernel Hilbert spaces and
Agler-Pick interpolation has been, starting with the work of Agler
\cite{MR1207393}, generalized from the classical setting of
$H^\infty(\mathbb D)$ ($\mathbb D$ the unit disk), to many
other algebras of functions.

In this paper we pursue realization formul{\ae} and Agler-Pick
interpolation in two directions.  First we consider an algebra of
functions on a semigroupoid $G$.  The precise definition of a
semigroupoid is given below.  In essence it can be thought of as an
ordered unital semigroup, though perhaps with more than one unit.  For
now the salient point is that the semigroupoid structure means that
the algebra product generalizes both the pointwise and convolution
products.  This setting has the advantage of being fairly concrete and
amenable to study using reproducing kernel Hilbert space ideas and
techniques while at the same time connecting with the theory of graph
$C^*$-algebras.

Secondly, we view the norm on the algebra as being determined by a
(possibly infinite) collection $\Psi$ of functions on $G$, referred to
as test functions.  Results on Agler-Pick interpolation (in the
classical sense) for both finite and infinite collections of test
functions with varying amounts of additional imposed structure can be
found in \cite{MR2106336,MR1882259} and this point of view
goes back at least to \cite{Ag-unpublished}.  A collection of test
functions determines a family of kernels, and vice versa.  This
duality between test functions and kernels will have a familiar feel
to those acquainted with Agler's model theory \cite{MR976842}.  The
advantage of such an approach is that it allows us to consider
interpolation problems on, for example, polydisks and multiply
connected domains \cite{DMnext}.

We should mention that Kribs and Power \cite{MR2076898,MR2063759}
introduce a somewhat more restrictive notion of a semigroupoid
algebra.  These are related to so-called quiver algebras of Muhly
\cite{MR1462686}, and are the nonselfadjoint analogues of the higher
rank graph algebras of Kumjian and Pask \cite{MR1745529}.  In these
papers order is either imposed through the presence of a functor from
the semigroupoid to ${\mathbb N}^d$, or by the assumption of freeness.
In either case, the resulting object is cancellative, and there is a
representation (related to our Toeplitz representation on
characteristic functions $\chi_a$; see
section~\ref{subsec:intro-standardrep}) in terms of partial isometries
and projections on a generalized Fock space with orthonormal basis
labelled by the elements of the semigroupoid.  The algebras of
interest in these papers are obtained as the weak operator topology
closure of the algebras coming from the left regular representation
(i.e., the projections and partial isometries mentioned above), and so
in a natural sense are the multiplier algebras for these Fock spaces.

The Kribs and Power semigroupoid algebras include the noncommutative
Toeplitz algebras first introduced in \cite{MR1129595}.  Pick and
Carath\'eodory interpolation has been considered in this context by
Arias and Popescu \cite{MR1749679} and Davidson and Pitts
\cite{MR1627901} (with some earlier work by Popescu on these and
related interpolation problems to be found in
\citep{MR1652931,MR1348353,MR1033608,MR1026074}), and somewhat more
generally by Jury and Kribs \cite{MR2153149}.  See also
\cite{MR2182589}.  In fact, while the commutant lifting theorem
unifies the classical Pick and Carath\'eodory-Fej\'er interpolation
problems, to our knowledge, Jury's PhD dissertation \cite{Jury-diss}
was the first to do so in terms of the positivity of kernels, and also
the first to give a concrete realization formula for the case of the
semigroup $\mathbb N$.  Recently, realization formul{\ae} in a
noncommutative setting have also been investigated in
\cite{MR2232307}.  Muhly and Solel \cite{MR2129649} have considered
Nevanlinna-Pick interpolation from the vantage of what they call Hardy
algebras, covering many of the examples mentioned above along with the
statement of a realization formula.

Interpolation problems on domains other than the unit disk in $\mathbb
C$ have been of long-standing interest.  On multiply connected
domains, the seminal work is that of Abrahamse \cite{MR532320}, with
further contributions to be found in
\citep{MR1386327,MR1367626,MR1707931,MR1818066,MR1691394,MR2163865}.
Regarding domains in $\mathbb C^n$, the fundamental paper of Agler
\cite{MR1207393} provides the foundation upon which most subsequent
work is based.  A sampling of papers of particular interest in this
direction includes \citep{MR1665697,MR1882259,MR1885440,MR2106336,
MR1722812,MR1880832,MR1748285,MR1846055,MR2050101,MR1794817,MR2069781,
MR2084918}.

In this paper we have for clarity restricted our attention to scalar
valued interpolation (although we stray a bit in the examples in
Section~\ref{sec:examples}).  We do not anticipate that the
generalization to the matrix case will provide any obstacles which
cannot be overcome with what are by now standard techniques.  Indeed
we have ensured that none of the proofs found below depend on the
commutativity of the coefficients of our functions, and it appears
likely that the results will continue to hold when the coefficients
come from, say, a norm closed subalgebra of a $C^*$-algebra.  This is
left for later work.

A few words about the organization of the paper.  The rest of
Section~1 outlines the basic tools used throughout: semigroupoids,
$\st$-products, Toeplitz representations, test functions and
reproducing kernels, the $C^*$-algebra generated by evaluations on the
set of test functions and its dual, and transfer functions.  This is
followed by a statement of the main results, which are the realization
and interpolation theorems.

In Section~2 we more closely study the $\st$-product, especially with
regards to inverses and positivity.

Section~3 begins with a consideration of the semigroupoid algebra
analogue of the Szeg\H{o} kernel, and highlights the close connection
between positivity of these kernels and complete positivity of the
$\st$-product map (a generalization of the Schur product map).  As
noted earlier in the introduction, multiplier algebras arising from a
single reproducing kernel are too restrictive for us, so we detail how
we will handle families of kernels and the associated families of test
functions.  Cyclic representations of the space of functions over
certain finite sets (they should be ``lower'' with respect to the
order on the semigroupoid) which are contractive on test functions are
shown to be connected to reproducing kernels.  This plays a crucial
role in the Hahn-Banach separation argument in the realization
theorem.

Given a positive object, an analyst's immediate inclination is to
factor.  The fourth section is devoted to a factorization result for
positive kernels on the dual of the $C^*$-algebra from Section~1, as
well as making connections to representations of this algebra.

Two other key items needed in the proof of the realization theorem
are taken up in Section~5.  The first is the cone of matrices
$\mathcal C_F$.  For the separation argument in the proof of the
realization theorem to work, we must know that $\mathcal C_F$ is
closed and has nonempty interior.  Closedness requires a surprisingly
delicate argument, and so occupies the bulk of the section.  We also
show that certain sets of kernels in the dual of the $C^*$-algebra
mentioned above are compact.

Sections 6 and 7 comprise the proof of the realization theorem and the
interpolation theorem.  The first implication of the proof of the
realization theorem is essentially the most involved part, but due to
all of the preparatory work in Sections~3--5, is dispensed with
quickly.  Other parts involve variations on themes which will be
familiar to those acquainted with recent proofs of interpolation
results.  These include an application of Kurosh's theorem, a lurking
isometry argument, and a fair amount of tedious calculation.  After
the proof of the realization theorem, the proof of the interpolation
theorem is almost an afterthought.

In Section~8 we turn briefly to a menagerie of examples, both old and
new.  Though we mention it in passing, we have postponed the
application to Agler-Pick interpolation on an annulus to a separate
paper for two reasons.  First, the argument is fairly long and
involves ideas and techniques unrelated to the rest of this paper; and
second, the underlying semigroupoid structure is that of Pick
semigroupoid (which is essentially trivial) and as such the version of
Theorem~\ref{thm:main} which is needed does not require the
semigroupoid overhead.  In any case, this section barely scratches the
surface of what is possible!

We would like to thank Robert Archer for his careful reading of the
paper, and the many useful comments and questions which have without a
doubt improved it.

\subsection{Semigroupoids}
\label{subsec:intro-semigroups}

There is no standard name in the literature for the sort of object on
which we want to define our function algebras.  The names ``small
category'' and ``semigroupoid'' are two commonly used terms, though
our definition differs somewhat from that standardly given for either
of these.  We have opted for the latter.

The term ``semigroupoid'' was originally coined by Vagner, as far as
we are aware \cite{MR0220856}.  Similar notions are familiar from the
theory of inverse semigroups (see for example, \cite{MR1694900} or
\cite{MR1724106}), and have been explored in connection with the
classification theory of $C^*$-algebras.  The use of semigroupoids in
the study of nonselfadjoint algebras originates with Kribs and Power
\cite{MR2076898}, though again, their use of the terminology is a bit
different from ours.

So let $G$ be a set with a function $X\subset G\times G \to G$, called
a \textit{partial multiplication} and written $xy$ for $(x,y)\in X$.
We define \textit{idempotents} as those elements $e$ of $G$ such that
$ex=x$ whenever $ex$ is defined and $ye=y$ whenever $ye$ is defined.
Note that these are commonly referred to as \textit{identities} in the
groupoid literature.

The following laws are assumed to hold:
\begin{enumerate}
\item (\textbf{associative law}) If either $(ab)c$ or $a(bc)$ is
  defined, then so is the other and they are equal.  Also if $ab$,
  $bc$ are defined, then so is $(ab)c$.
\item (\textbf{existence of idempotents}) For each $a\in G$, there
  exist $e,f\in G$ with $ea=a=af$.  Furthermore if $e\in G$ satisfies
  $e^2=e$, then $e$ is idempotent.
\item (\textbf{nonexistence of inverses}) If $a,b\in G$ and $ab=e$
  where $e$ is idempotent, then $a=b=e$.
\item (\textbf{strong artinian law}) For any $a\in G$ the cardinality
  of the set $\{z,b,w : zbw=a\}$ is finite.  Moreover there is an $N<
  \infty$ such that $\sup_{a,c\in G} \textrm{card} \{b\in G : cb = a
  \} \leq N$.
\end{enumerate}

Hereafter we refer to a set $G$ with a partially defined
multiplication with all of the properties so far listed as a
\textit{semigroupoid}.

Since we have associativity, we can mostly forget parentheses.  If we
were to reverse the third law (so that every element has an inverse),
then the first three rules would comprise the definition of a
groupoid.  The strong artinian law is related to the (partial) order
which we eventually impose on our semigroupoid.  The first part of it
ensures that the multiplication we will define for functions over the
semigroupoid is well defined, while the second part guarantees the
existence of at least one collection of test functions, or
equivalently, that the associated collection of reproducing kernels is
nontrivial.  It does so by restricting how badly non-cancellative the
semigroupoid can be.  Alternately, the strong artinian law could be
replaced by the condition that for each $a\in G$ the set
$\{z,b,w:zbw=a\}$ is finite and a hypothesis about the existence of a
collection of test functions (see
Section~\ref{subsec:intro-test-functions}).

There is one other rule which it is useful to state, though it follows
from those already given:
\begin{itemize}
\item[(5)] (\textbf{strong idempotent law}) If $zaw = a$, then $z$
  and $w$ are idempotents.
\end{itemize}
To see that this is a consequence of our other laws, first note that
$zaw=a$ means that $z^naw^n=a$ for $n\in \mathbb N$.  The strong
artinian law implies that only finitely many of the $z^n$ are
distinct.  In particular, there is an $M > 0$ such that $(z^{2^M})^2 =
z^{2^j}$ for some $j\leq M$.  Let $h_1=z^{2^j}$, $h_2 = h_1^2$, and so
on, with $h_m = h_{m-1}^2 = z^{2^M}$ (so $h_1 = h_m^2$).  Clearly
$h_jh_k=h_kh_j$ for all $j,k$.  Hence
\begin{equation*}
  (h_1h_2\cdots h_m)^2 = h_1^2 \cdots h_m^2 = h_2 \cdots h_m h_1 =
  h_1h_2\cdots h_m,
\end{equation*}
and so $h_1h_2\cdots h_m = z^{2^j+\cdots + 2^M}$ is idempotent.  Since
there are no inverses, this implies that $z$ is idempotent.  Likewise
$w$ is idempotent.  Note that (5) implies our assumption that $e$ is
idempotent if $e^2=e$.

If $ea=a$ then $e$ is unique, since if $e'a=a$, then $a = ea = e(e'a)
= (ee')a$, implying that $ee'$ is defined.  But then since $e$ and
$e'$ are assumed to be idempotents, $e=ee'=e'$.

From the definition we have $e^2=e$ means that $e$ is idempotent.  On
the other hand, if $a=ea$ then $a = e(ea) = e^2a$, and so $e^2$ is
defined, and by uniqueness, $e^2=e$.  Also if $e$ and $f$ are
idempotents and $ef$ is defined, then $e=ef=f$.

The product $ab$ exists if and only if there is an idempotent $f$ such
that $af$, $fb$ are defined.  For if such an $f$ exists, then by
associativity, $(af)b = a(fb) = ab$, while conversely, if $ab$ is
defined, then there is an idempotent $f$ such that $af=a$ and so
$(af)b=a(fb)$ is defined and so $fb$ is defined.

Based on these observations, it is common to view a set with a partial
multiplication verifying the first two rules as a sort of generalized
directed graph with the vertices representing the idempotents, though
because we have not assumed any cancellation properties, this analogy
is imperfect.

We define subsemigroupoids in the obvious way.  In particular, a
subset $H$ of a semigroupoid $G$ will be a subsemigroupoid if whenever
$a,b\in H$ and $ab$ makes sense in $G$ then $ab\in H$, and for all
$a\in H$ the idempotents $e,f$ such that $ea=a=af$ are also in $H$.

We put a partial order on a semigroupoid $G$ as follows: say that $b
\leq a$ if there exist $z,w\in G$ such that $a=zbw$.  By the existence
of idempotents, $a\leq a$.  Transitivity is likewise readily verified.
If $a\leq b$ and $b\leq a$ then $a=zbw$, $b=z'aw'$ and so
$a=(zz')a(w'w)$.  Then by the strong idempotent law $zz'$ and $w'w$
are idempotent.  But then by the nonexistence of inverses, $z$, $z'$,
$w$ and $w'$ are idempotents and so $a=b$.  Other partial orders are
considered in Section~\ref{sec:more-order-semigr}.

By this definition, and the existence of idempotents,
if $a=bc$, then both $b$ and $c$ are less than or equal to $a$.  Also,
by the nonexistence of inverses and uniqueness of idempotents, the
idempotents comprise the minimal elements of $G$.  We write $G_e$ for
the collection of idempotents.

We say that a set $F\subset G$ is \textit{lower} if $a\in F$ and
$b\leq a$ then $b\in F$.  Observe that for a lower set $F$, $F_e =
F\cap G_e \neq \emptyset$.  Note too that if $H$ is a finite subset of
$G$, then there is a finite lower set $F\supset H$: simply let $F=\{a:
\text{ there exists a } b\in H \text{ such that } a\le b\}$.

\subsubsection{Examples}
\label{subsubsec:examples}

We list here several important examples of semigroupoids.

\begin{enumerate}
\item Let $G$ be a set, and assume $G_e = G$ (so all elements are
  idempotent).  We refer to such semigroupoids as \textit{Pick
    semigroupoids}.
\item Let $G = \mathbb{N} =0,1,2,\ldots$ with the product $ab = a+b$.
  $G$ is in fact a commutative cancellative semigroup with idempotent
  $0$.
\item The last example obviously generalizes to $\mathfrak{F}_n$, the
  free (noncommutative) monoid on $n$ generators.  This in turn is a
  special case of what we term the \textit{Kribs-Power semigroupoids}
  \cite{MR2076898}, which are defined as follows.  Let $\Lambda$ be a
  countable directed graph.  The semigroupoid $\mathfrak F^+(\Lambda)$
  determined by $\Lambda$ comprises the vertices of $\Lambda$, which
  act as idempotents, and all allowable finite paths in $\Lambda$,
  with the natural concatenation of allowable paths in $\Lambda$
  defining the partial multiplication.  In particular, $\mathfrak
  F^+(\Lambda) = \mathfrak{F}_n$ when $\Lambda$ is a directed graph
  with one vertex and $n$ distinct loops.
\end{enumerate}

\subsection{The convolution products}
\label{subsec:intro-products}
 
The product on $G$ naturally leads to a product on functions over
lower sets $F\subset G$ in one or more variables.

\subsubsection{The $\st$-product for functions of one variable}
\label{subsubsec:intro-starproduct}

Let $F$ be a lower subset of $G$.  There is a natural algebra
structure on the set $P(F)$ of functions $f:F\to\mathbb C$ which we
call the semigroupoid algebra of $F$ over $\mathbb C$.  Addition of
$f,g\in P(F)$ is the usual pointwise addition of functions and the
product is defined by
\begin{equation*}
  (f\s g)(a) = \sum_{rs=a} f(r)g(s),
\end{equation*}
which makes sense because of the artinian hypothesis on $G$ and
the assumption that $F$ is lower.

The multiplicative unit of $P(F)$ is given by
\begin{equation*}
  \delta(x) = 
  \begin{cases}
    1 & x\in F_e,\\
    0 & \text{otherwise.}
  \end{cases}
\end{equation*}
The distributive and associative properties are readily checked, so we
have an algebra.  A function $f$ is invertible if and only if $f(x)$
is invertible for all $x\in F_e$.  The proof follows the same lines as
in the matrix case given below, so we do not give it here.

If $a\in F^\prime\subset F$ and $F^\prime$ is itself lower, then
\begin{equation*}
  ( f_{|F^\prime} \s g_{|F^\prime})(a)= (f\s g)(a).
\end{equation*}
Hence, we can be lax in specifying our lower set and usually act as if
it is finite.

Later we have need for powers of functions with respect to the
$\st$-product.  To avoid confusion, for a function $\varphi$ on $G$, we
let $\varphi^{n_\st}$ denote the $n$-fold $\st$-product of $\varphi$
with itself.

As it happens, it is unimportant that a function over $F$ map into
$\mathbb C$.  For instance, the $\st$-product clearly generalizes to
functions $f,g:F\to \mathfrak C$, where $\mathfrak C$ is a
$C^*$-algebra.

There will be times when we will want to interchange $r$ and $s$ in
the definition of the convolution product.  Over $\mathbb C$ or, more
generally, any commutative $C^*$-algebra $\mathfrak C$ this simply
changes $f\s g$ into $g\s f$.  But in the noncommutative case this
will not work.  Hence we introduce the notation
\begin{equation*}
  (f\sh g)(a) = \sum_{rs=a} f(s)g(r).
\end{equation*}
For the $\sth$-product the multiplicative unit remains $\delta$, the
associative and distributive laws continue to hold, and $f$ is
invertible with respect to this product if and only if $f(x)$ is
invertible for all $x\in F_e$.  We write $f^{-1_{\st}}$ and
$f^{-1_{\hat{\star}}}$ for the $\st$-inverse and $\sth$-inverse of
$f$, respectively.  By considering $f^{-1_{\hat{\star}}} \sh f \s
f^{-1_{\st}}$, we see that $f^{-1_{\hat{\star}}} = f^{-1_{\st}}$.

Another useful and easily checked property relating the two products
is that
\begin{equation}
  \label{eq:5}
  (f\sh g)^* = g^*\s f^*.
\end{equation}
Consequently $(f\sh g)(g^*\s f^*) \geq 0$.  We also have
$(f^{-1_{\st}})^* = f^{-1_{\hat{\star}}}$.

In the examples listed above, the $\st$-product is just pointwise
multiplication for Pick semigroupoids.  For the second example, it is
the usual convolution.

\subsubsection{The $\st$-product for matrices}
\label{subsubsec:intro-defin-prop}

The following bivariate version of the convolution product is the
canonical generalization of Jury's product \cite{Jury-diss} to
semigroupoids.

For a lower set $F$, let $M(F)$ denote the set of functions $A:F\times
F\to \mathbb C$.  When $F$ is finite, thinking of elements of $M(F)$
as matrices (indexed by $F$), the notation $A_{a,b}$ is used
interchangeably with $A(a,b)$.  The set of functions from $F\times F$
to $X$ will be denoted $M(F,X)$.

\begin{definition}
  Let $F$ be a lower set and suppose  $A, B\in M(F)$.
  Define $A\s B$ by
  \begin{equation*}
    (A\s B)(a,b) = \sum_{pq=a} \sum_{rs=b} A(p,r)B(q,s).
  \end{equation*}
\end{definition}

Once again, the artinian hypothesis on $G$ guarantees the product is
defined.  Further, $(A\s B)(a,b)$ does not actually depend upon the
lower set $F$ which contains $a$ and $b$.  In particular, since there
is always a finite lower set containing $a$ and $b$ (just take the
union of the set of elements less than or equal to $a$ and those less
than or equal to $b$), this product can and will be interpreted as a
$\st$-product of matrices.

The assumption that the entries of $A$ and $B$ are in $\mathbb C$ is
not important, and we will at times use the $\st$-product when the
entries are in other algebras.  The $\st$ notation should cause no
confusion, since in essence the $\st$-product is the bivariate analogue
of the convolution product.  Indeed, it is clear that the $\st$-product
could be defined for functions of three or more variables as well,
though we have no need for this here.

Unlike Jury's $\st$-product, ours will not necessarily be commutative
(though this will be the case if $G$ is commutative).  In the special
example of the Pick semigroupoids, the $\st$-product is just the matrix
Schur product.

As with functions we can also define the $\sth$-product of matrices:
  \begin{equation*}
    (A\sh B)(a,b) = \sum_{pq=a} \sum_{rs=b} A(q,s)B(p,r).
  \end{equation*}
  Over $\mathbb C$ and any other commutative algebra, $A\sh B = B\s
  A$.  However we will need both products in a noncommutative setting.

Define $[1]$ (or $[1]_F$ if we want to make absolutely clear the
lower set involved) to be the matrix defined by $[1]_{a,b}$ is $1$ for
$a,b$ both elements of $F_e$ and zero otherwise.  It
is easy to see that for $A\in M(F)$, $A\s [1] = [1]\s A =
A$.  Note too that we can factor $[1] = \delta \delta^*$.

\subsection{The Toeplitz representation}
\label{subsec:intro-standardrep}

Let $\varphi$ be a function on a (finite) lower set $F$.  Define the
associated \textit{(left) Toeplitz representation $\sT (\varphi)$} by
\begin{equation*}
  [\sT (\varphi)]_{a,b} = 
  \begin{cases}
    \sum_c\varphi(c), & cb=a; \\
    0 & \text{otherwise}.
  \end{cases}
\end{equation*}
We drop the ``left'' hereafter, though we could also consider a right
Toeplitz representation with $bc=a$ rather than $cb=a$.  It is
expedient to assume $F$ is finite to ensure that $\sT (\varphi)$ is
bounded, though the definition makes sense formally when $F$ is not
finite and in many interesting cases yields a bounded operator.

As defined, $\sT(\varphi)$ is a mapping of $F\times F$ into $\mathbb
C$.  Let $\mathbb C^F$ denote the column vector space (of dimension
equal to the cardinality of $F$) with positions labelled by the
entries of $F$, which is naturally isomorphic to $P(F)$ as a vector
space.  Then for a function $f\in P(F)$, viewed as an element
of $\mathbb C^F$,
\begin{equation}
  \label{eq:2}
  (\sT (\varphi) f)(a) = \sum_b [\sT (\varphi)]_{a,b} f (b) =
  \sum_{cb=a} \varphi(c) f(b) = (\varphi \star f)(a).
\end{equation}
In this way $\sT(\varphi)$ is an operator in $B(\mathbb C^F)$ and the
mapping $P(F)\ni \varphi \mapsto \sT(\varphi)\in B(\mathbb C^F)$ is a
representation.  Indeed, it essentially acts as the left regular
representation of $P(F)$.  The use of the notation $M(F)$ to denote
either $B(\mathbb C^F)$ or functions from $F\times F\to \mathbb C$
should be clear from the context.  Further, since the Toeplitz
representation depends upon the lower set $F$ in a consistent way, it
should cause no serious harm that the notation $\sT (\varphi)$ makes
no reference to $F$.  On occasions when we need to make the dependence
on $F$ explicit, we will write $\sT^F$ for $\sT$.

In the case that $G$ is the semigroupoid $\mathbb N$, $\sT(\varphi)$
is precisely the Toeplitz matrix associated with the sequence
$\{\varphi(j)\}$.  At the other extreme, when $G$ is a Pick
semigroupoid, $\sT(\varphi)$ is the diagonal matrix with diagonal
entries $\varphi(a)$ for $a\in G$ which, despite our terminology,
seems very un-Toeplitz like!

When $G$ is the Kribs-Power semigroupoid $\mathfrak F^+(\Lambda)$
determined by a countable directed graph $\Lambda$, a lower set $F$ of
$\mathfrak F^+(\Lambda)$ is closed under taking left and right
subpaths.  The vector space $P(F)$ may be regarded as a subspace of
the generalized Fock space $\mathcal H_\Lambda$ over $\Lambda$ (see
\cite{MR2153149}).  In this interpretation $\sT^F(\chi_w)$ is the
compression to $P(F)$ of the partial creation operators indexed by $w
\in F$.  More generally, $\{\chi_v~:~v \in \mathfrak F^+(\Lambda)\}$
may be thought as an orthonormal basis of $\mathcal H_\Lambda$ and
$\sT$ behaves as a representation $\mathfrak F^+(\Lambda) \to
B(\mathcal H_ \Lambda)$.  Hence, the weak operator topology closed
subalgebra generated by the family $\{\sT(\chi_v)\}$ is the free
semigroupoid algebra of Kribs and Power \cite{MR2076898}, which
includes, as particular case, the noncommutative Toeplitz algebra
\cite{MR1627901,MR1749679,MR1675377}.  The set of
generators can be restricted, as we see later.

Even when $F$ is not necessarily finite, $\sT$ still behaves
formally as a representation, but of course it need not
be the case that $\sT(\varphi)$ is bounded.

It is also possible to work with the $\sth$-product.  Presumably, there
is a distinction between a collection of test functions, defined
below with respect to the $\st$-product, and those with respect to the
$\sth$-product, though we do not develop this.

\subsection{Test functions}
\label{subsec:intro-test-functions}

For a function $\varphi$ on $G$, recall that $\varphi^{n_\st}$ denotes
the $n$-fold $\st$-product of $\varphi$ with itself, $n=1,2,\ldots$.

\begin{definition} \rm
  \label{def:test-functions}
  A collection $\Psi$ of functions on $G$ into $\mathbb C$ is a
  \textit{collection of test functions} if
  \begin{enumerate}
  \item [(i)] for each finite lower set $F\subset G$ and $\psi\in\Psi$,
    $\|\sT(\psi)\|\le 1$;
  \item [(ii)] for each $a\in G_e$, 
    \begin{equation*}
      \lim_{n\to\infty} \psi^{n_\st}(a) = \lim_{n\to\infty} \psi^n(a) =
      0,
    \end{equation*}
    uniformly in $\psi$; and
  \item [(iii)] for each finite lower set $F$, the algebra generated by
    $\PrF=\{\psi|_F: \psi\in\Psi\}$ is all of $P(F)$.
  \end{enumerate}
\end{definition}

The condition $\PrF$ generates $P(F)$ is not essential.  It does
however simplify statements of results.

Given $x\in G$, let $f$ be the unique idempotent so that $xf=x$.
Since,
\begin{equation*}
  [\sT(\psi)]_{x,f} = \sum_{cf=x}\psi(c)= \psi(x),
\end{equation*}
item (i) says that $|\psi(x)|\le 1$ for each $x\in G$.  By the same
reasoning, if $\psi_1,\psi_2\in\Psi$, then
\begin{equation*}
  |\psi_1(x)-\psi_2(x)|\le \|\sT(\psi_1)-\sT(\psi_2)\|.
\end{equation*}

Item (ii) says that for each $a\in G_e$ and $\epsilon>0$ there is an
$N$ so that for all $n\ge N$ and $\psi\in \Psi$, $|\psi^n (a)| <
\epsilon$, and so for fixed $a\in G_e$, $\sup_{\psi\in\Psi} |\psi (a)|
< 1$.  Furthermore, for any $a\in G$, we automatically obtain
$\lim_{n\to\infty} \psi^{n_\st}(a) =0$.  This follows from a
straightforward counting argument estimating the maximum number of
ways of writing $a\in G$ as a product of $n$ elements.  Assume
$\mathrm{card} \{ b : b \leq a\} = r$ (which is finite since $G$ is
artinian) and that $n \gg r$.  Let $c = \max_{b\leq a} |\psi(b)|$, and
$c_e$ be the maximum of $|\psi(b)|$ over all idempotents less than or
equal to $a$.  As noted above, $c_e < 1$.  In the product of $n$
terms, there are at most $\binom{n}{r}$ ways of choosing which of the
at most $r$ terms are not idempotent, and then at most $r^r$ ways of
choosing these terms.  The nonidempotents act as separators between at
most $r+1$ blocks of idempotents.  Within each block of idempotents,
each term must be the same idempotent (since the product of unequal
idempotents is not defined).  So there are at most $r^{r+1}$ ways of
choosing which idempotent is in each block.  Consequently
\begin{equation*}
  |\psi^{n_\st}(a)| \leq \binom{n}{r}r^r c^r r^{r+1} c_e^{n-r} \leq
  r^{2r+1} (c/c_e)^r n^r c_e^n ,
\end{equation*}
which clearly goes to zero as $n\to\infty$.

For a given semigroupoid $G$ it is legitimate to wonder if there
actually exists any family of test functions.  It so happens that the
strong artinian condition in the definition of a semigroupoid ensures
this.  Let $\kappa = \sup_{a,c\in G} \textrm{card} \{ b\in G : cb = a
\}$, which we have assumed is finite.  Let $\psi_0:G\to \mathbb D$
with $\psi_0|_{G_e}$ injective and $\psi_0|_{G\backslash G_e} = 0$.
(This assumes the cardinality of $G_e$ is less than or equal to that
of the continuum --- it is only slightly more trouble to handle the
more general case.)  Let $\Psi_s = \{ \frac{1}{\kappa} \chi_c : c\in
G\backslash G_e\} \cup \{\psi_0\}$.  Then $\Psi_s$ can be shown to be
a collection of test functions (here $\chi_c(x)$ equals $1$ if $x=c$
and zero otherwise).  In particular, the condition $\kappa < \infty$
for all $c\in G$ will hold if $G$ is right cancellative (so in
particular, for Kribs-Power semigroupoids).  In the case $G=G_e$, this
choice of test functions will ultimately correspond to $B(G)$, the
normed algebra of all bounded functions on $G$.

\subsection{Test functions and reproducing kernel Hilbert spaces}
\label{subsec:intro-rkhs}

Let $F\subset G$ be a lower set.  A function $\mathbf{k}:F\times F\to
\mathbb C$ is a positive kernel if for each finite subset $A\subset F$
the matrix $[\mathbf{k}(a,b)]_{a,b\in A}$ is positive (i.e., positive
semidefinite).
 
More generally, it makes sense to speak of a kernel with values in the
dual of a $C^*$ algebra.  If $\mathfrak B$ is a $C^*$-algebra with
Banach space dual $\mathfrak B^*$, then a function $\Gamma :F\times
F\to \mathfrak B^*$ is positive if for each finite subset $A\subset F$
and each function $f:A\to \mathfrak B$,
\begin{equation*}
   \sum_{x,y\in A} \Gamma(x,y)(f(x)f(y)^*) \ge 0.
\end{equation*}
In the sequel, unless indicated otherwise, kernels take their values
in $\mathbb C$.

Given a set of test functions $\Psi$ let ${\mathcal K}_\Psi$ denote
the collection of positive (i.e., positive semidefinite) kernels
$\mathbf{k}$ on $G$ such that for each $\psi\in\Psi$, the kernel
\begin{equation}
  \label{eq:defkpsi}
  G\times G \ni (x,y)\mapsto   \mathbf{k}_\psi(x,y) = \left(([1]-\psi\psi^*)\s
    \mathbf{k}\right)(x,y)
\end{equation}
is positive.  Here $[1]-\psi \psi^*$ is the function defined on
$G\times G$ by $([1]-\psi\psi^*)(p,q)=[1]_{p,q}-\psi(p)\psi(q)^*$ so
that the right hand side of equation \eqref{eq:defkpsi} is the
$\st$-product of the functions (or matrices indexed by $G$) $[1]-\psi
\psi^*$ and $\mathbf{k}$, evaluated at $(x,y)\in G\times G$.

The set ${\mathcal K}_\Psi$ is nonempty, since it at least contains
$\mathbf{k}=0$.  More importantly, from the hypothesis that $\Psi$ is
a family of test functions and the strong artinian law, it also
contains the kernel $s:G\times G\to \mathbb C$ given by $s(x,y)=1$ if
$x=y$ and $0$ otherwise, which is strictly positive definite.  We call
$s$ the \textit{Toeplitz kernel}.

Let us verify that $s\in {\mathcal K}_{\Psi_s}$ for the collection of
test functions $\Psi_s$ constructed in the last subsection.  For the
test function $\psi_0$,
\begin{equation*}
  \begin{split}
    (\psi_0\psi_0^* \s s)(a,b) & = \sum_{\scriptstyle
      pq=a\atop\scriptstyle rt=b} \psi_0(p)\psi_0^* (r) s_{q,t} =
    \sum_{\scriptstyle p,r\in G_e\atop{\scriptstyle
        pq=a\atop\scriptstyle rt=b}} \psi_0(p)\psi_0^* (r) s_{q,t}\\
    & = \begin{cases}
      \psi_0(p)\psi_0^* (p), & p\in G_e,\ a=b,\ pa=a\\
      0 & \text{otherwise.}
    \end{cases}    
  \end{split}
\end{equation*}
Hence since $[1]\s s = s$, $([1]-\psi_0\psi^*_0)\s s$ is a diagonal
matrix with entries of the form $1 - \psi_0(p)\psi^*_0(p) \geq 0$, and
so is positive.  On the other hand suppose $\psi_c =
\frac{1}{\kappa} \chi_c$, $c\in G\backslash G_e$.  Then
\begin{equation*}
  (\psi_c\psi_c^* \s s)(a,b) = \sum_{\scriptstyle
    pq=a\atop\scriptstyle rt=b} \psi_c(p)\psi_c^* (r) s_{q,t} =
  \tfrac{1}{\kappa^2}{\left(\sum_{cq=a=b} 1\right)}^2 \in [0,1],
\end{equation*}
and so $([1]-\psi_c\psi^*_c)\s s$ is a positive diagonal matrix.

The kernels determined by a family of test functions $\Psi$ in turn
give rise to a normed algebra of functions on $G$.  Let $\HP$ denote
those functions $\varphi:G\to \mathbb C$ such that there exists a
$C>0$ such that for each $\mathbf{k}\in \mathcal K_\Psi$, the kernel
\begin{equation*}
  G\times G \ni (x,y) \mapsto ((C^2[1]-\varphi\varphi^*)\s \mathbf{k})(x,y)
\end{equation*}
is positive.  The infimum of all such $C$ is the norm of $\varphi$.
With this norm $\HP$ is a Banach algebra under the convolution
product.  By construction $\Psi$ is a subset of the unit ball of
$\HP$.

There is a duality between kernels and test functions in Agler's model
theory \cite{Ag-unpublished,MR976842}.  Roughly, the idea
is, given a collection $\mathcal K$ of positive kernels on $G$, to let
$\Psi=\mathcal K^\perp$ denote those functions $\psi\in G$ such that
for each $\mathbf{k}\in\mathcal K$, the kernel
\begin{equation*}
   G\times G\ni (x,y)\mapsto (([1]-\psi \psi^*)\s \mathbf{k})(x,y)
\end{equation*}
is positive.  In the case that Agler considers, where the semigroupoid
consists solely of idempotents (i.e., a Pick semigroupoid), mild
additional hypotheses on $\mathcal K$ guarantee that $\Psi$ is a
family of test functions, in which case $\mathcal K_\Psi =\mathcal
K^{\perp\perp}$.

\subsection{The evaluation $E$ and $C^*$-algebra $\CT$}
\label{subsec:intro-evaluation}

Let $\Psi$ be a given collection of test functions and $C_b(\Psi)$ the
continuous functions on $\Psi$, where $\Psi$ is compact in the bounded
pointwise topology.  Define $E\in B(G,C_b(\Psi))$ (the bounded
functions from $G$ to $C_b(\Psi)$) by
\begin{equation*}
  E(x)(\psi) = \psi(x), \qquad \psi\in \Psi,
\end{equation*}
with
\begin{equation*}
  \|E(x)\| = \sup_{\psi\in \TT} \{|E(x)(\psi)|\}.
\end{equation*}
So $E(x)$ is the evaluation map on $\TT$, $\|E(x)\|<1$ for each $x\in
G_e$ and $\|E(x)\| \leq 1$ otherwise.

Since evidently the collection $\{E(x):x\in G\}$ separates points and
we include the identity, the smallest unital $C^*$-algebra containing
all the $E(x)$ is $C_b(\Psi)$.  For convenience, we denote this
algebra as $\CT$.

\subsection{Colligations}
\label{subsec:colligations}

Following \cite{MR2106336} we define a $\CT$-unitary colligation
$\Sigma$ to be a triple $\Sigma= (U,\mathcal E,\rho)$ where $\mathcal
E$ is a Hilbert space,
\begin{equation*}
  U=\begin{pmatrix} A & B \\  C & D \end{pmatrix} : 
  \begin{matrix} \mathcal E \\ \oplus \\ \mathbb C \end{matrix} \to 
  \begin{matrix} \mathcal E \\ \oplus \\ \mathbb C \end{matrix}
\end{equation*}
is unitary, and $\rho:\CT\to B(\mathcal E)$ is a unital
$*$-representation.  The transfer function associated to $\Sigma$ is
\begin{equation} 
  \label{eq:transfer}
  W_\Sigma(x)= (D\delta + C(\rho (E) \s (\delta -A\rho(E))^{-1_\st}
  \s (B\delta))(x).
\end{equation}
Observe that this looks like the standard transfer function over a
Pick semigroupoid.

\subsection{The main event}
\label{subsec:intro-main}

We now state the realization theorem for elements of the unit ball of
$\HP$ and a concomitant interpolation theorem.
 
\begin{theorem} [\textbf{Realization}]
  \label{thm:main}
  If $\Psi$ is a collection of test functions for the semigroupoid
  $G$, then the following are equivalent:
  \begin{itemize}
  \item[(i)] $\varphi\in \HP$ and
    $\|\varphi\|_{\HP} \le 1$;
  \item[(iiF)] for each finite lower set $F\subset G$ there exists a
    positive kernel $\Gamma:F\times F\to \CTs$ so that for all $x,y\in
    F$
    \begin{equation*}
      ([1]-\varphi\varphi^*)(x,y)= (\Gamma\sh([1]-EE^*))(x,y);
    \end{equation*}
  \item[(iiG)] there exists a positive kernel $\Gamma:G\times G\to
    \CTs$ so that for all $x,y\in G$
    \begin{equation*}
      ([1]-\varphi\varphi^*)(x,y)= (\Gamma\sh([1]-EE^*))(x,y);
      \text{ and}
    \end{equation*}
  \item[(iii)] there is a $\CT$-unitary colligation $\Sigma$ so that
    $\varphi=W_\Sigma$.
  \end{itemize}
\end{theorem}

\begin{theorem}[\textbf{Agler-Jury-Pick Interpolation}]
  \label{thm:main-interpolate}
  Let $F$ be a finite lower set and suppose $f\in P(F)$.  The
  following are equivalent:
  \begin{enumerate}
  \item[(i)] There exists $\varphi\in \HP$ so that 
   $\|\varphi\|_{\HP} \le 1$ and     $\varphi|_F=f$;
  \item[(ii)] for each $\mathbf{k}\in\mathcal K_\Psi$, the kernel
    \begin{equation*}
      F\times F \ni (x,y)\mapsto (([1]-ff^*)\s \mathbf{k})(x,y)
    \end{equation*}
    is positive; 
  \item[(iii)] there is a positive kernel $\Gamma : F\times F \to
    \CT^*$ so that for all $x,y\in F$
    \begin{equation*}
      ([1]-ff^*)(x,y) =(\Gamma \sh ([1]-EE^*))(x,y).
    \end{equation*}
  \end{enumerate}
\end{theorem}

\begin{remark} \rm
  \label{rem:AJP-int-realization}
  Of course, item (i) of Theorem~\ref{thm:main-interpolate} combined
  with item (iii) of Theorem~\ref{thm:main} says that $\varphi$ in
  Theorem~\ref{thm:main-interpolate} has a $\CT$-unitary transfer
  function representation.
\end{remark}

The hypothesis that $\PrF$ generates all of $P(F)$ means that the
representation $\pi:\HP \to P(F)$ which sends $\varphi$ to
$\varphi|_F$ is onto and identifies $P(F)$ with the quotient
$\HP/\mbox{ker}(\pi)$.  Theorem \ref{thm:main-interpolate} can be
interpreted as identifying the quotient norm.



\section{Further properties of the $\st$-products}
\label{sec:star-product}

\subsection{The convolution products}
\label{sec:convolution-product}

The convolution products over finite lower sets $F$ can be related to
the tensor product of matrices as follows.  Take $V:\mathbb C^n \to
\mathbb C^n\otimes \mathbb C^n$, where $n=\mathrm{card}(F)$, such that
$Ve_a = \sum_{pq=a} e_p\otimes e_q$, $\{e_k\}$ the standard basis for
$\mathbb C^n$ labelled with the elements of $F$, and extending by
linearity.  Then $f\s g = V^* (f\otimes g)$ and $f\sh g = (g\otimes
f)V$.  Note that $V$ is an isometry only in the case that $F = F_e$,
in which case the convolution products become the pointwise product.
In all other cases it still has zero kernel and in fact maps
orthogonal basis vectors $e_a$ and $e_b$ to orthogonal vectors, though
it generally acts expansively.

\subsection{The matrix $\st$-product---basic properties and an
  alternate definition}
\label{sec:defin-prop}

Straightforward calculations show that the various associative and
distributive laws hold for the bivariate $\st$-product.  Here, for
example, is the proof that $C\s(A\s B) = (C\s A)\s B$:
\begin{eqnarray*}
  [C\s(A\s B)]_{\mu,\nu} &=& \sum_{lm=\mu}\sum_{jk=\nu}
  C_{l,j}[A\s B]_{m,k} \\
  &=& \sum_{lm=\mu}\sum_{jk=\nu} C_{l,j}\sum_{pq=m}\sum_{rs=k}
  A_{p,r}B_{q,s} \\
  &=& \sum_{\scriptstyle l,p,q\atop\scriptstyle l(pq)=\mu}
  \sum_{\scriptstyle j,r,s\atop\scriptstyle j(rs)=\nu}
  C_{l,j}A_{p,r}B_{q,s},
\end{eqnarray*}
while
\begin{eqnarray*}
  [(C\s A)\s B)]_{\mu,\nu} &=& \sum_{iq=\mu}\sum_{ns=\nu}
  [C\s A]_{i,n}B_{q,s} \\
  &=& \sum_{iq=\mu}\sum_{ns=\nu}\sum_{lp=i}\sum_{jr=n} 
      C_{l,j}A_{p,r}B_{q,s} \\
  &=& \sum_{\scriptstyle l,p,q\atop\scriptstyle (lp)q=\mu}
  \sum_{\scriptstyle j,r,s\atop\scriptstyle j(rs)=\nu}
  C_{l,j}A_{p,r}B_{q,s},
\end{eqnarray*}
and since $l(pq) = (lp)q$, $j(rs) = (jr)s$, the two are equal.

There is an alternate equivalent definition of the $\st$-product, just
as with the convolution products.  Take $V$ defined as in the last
subsection.  Then it is easy to check that
\begin{equation}
  \label{eq:6}
  A\s B = V^*(A\otimes B)V.
\end{equation}
The Schur product is the matrix analogue of the pointwise product of
functions in which case $V$ is isometric, though otherwise it will not
be.  From this formulation it is clear that the $\st$-product is
continuous.

Another important property which the $\st$-product shares with the
Schur product is that if $A,B \in M(F)$ are positive, then so is $A\s
B$.  This follows immediately from the fact that $A\otimes B \geq 0$
if $A, B \geq 0$.  Similarly, since the tensor product of selfadjoint
matrices is selfadjoint, the $\st$-product of selfadjoint matrices is
selfadjoint.

\subsection{Positivity and the $\st$-product}
\label{subsec:positivity}

It should be emphasized that unlike with ordinary matrix
multiplication, the inverse with respect to the $\st$-product of a
positive matrix need not be positive.  This is already clear when
considering Schur products, but we illustrate with another simple
example.  Suppose that $e,a\in G$ with $e$ idempotent and $eae=a$.
Consider the matrix $A = \begin{pmatrix}1&0 \\0&c\end{pmatrix}$ where
$c>0$ and the first row and column is labelled by $e$ while the second
is labelled by $a$.  An easy calculation shows that $A^{-1_\st} =
\begin{pmatrix}1&0 \\0&-c\end{pmatrix}$.

The $\st$-product behaves somewhat unexpectedly with respect to
adjoints (at least if you forget its connection to the tensor
product).  Using the formulation of the $\st$-product given in
\eqref{eq:6}, we see that $(A\s B)^* = V^*(A\otimes B)^*V =
V^*(A^*\otimes B^*)V = A^*\s B^*$.  However with regard to inverses
and adjoints, $[1] = [1]^* = (A\s A^{-1_\st})^* = V^*(A^*\otimes
(A^{-1_\st})^*)V$, and so by uniqueness of the inverse, $A^*$ is
invertible if $A$ is and $(A^*)^{-1_\st} = (A^{-1_\st})^*$.
Consequently we see that if $A$ is selfadjoint and invertible, then
$A^{-1_\st}$ is selfadjoint.

Let $F$ be a finite lower set.  An $A\in M(F)$ gives rise to the
$\st$-product operator $S_A:M(F)\to M(F)$ given by $S_A(B)=A\s B =V^*
(A\otimes B)V$.  The argument in Paulsen's book (\cite{MR1976867},
Theorem~3.7) which shows that Schur product with a positive matrix
gives a completely positive map carries over with the obvious
modifications to show that $S_A$ is completely positive.  In
particular, the cb-norm of $S_A$ is given by $\|A\s 1\|$, where $1\in
M(F)$ is the identity (not the $\st$-product identity).

All of the above carries over in total to the $\sth$-product, with a
small change in the definition in terms of the tensor product, where
we have
\begin{equation*}
  A\sh B = V^*(B\otimes A)V.
\end{equation*}

\subsection{More on order on semigroupoids}
\label{sec:more-order-semigr}

The following lemmas give general properties of an artinian order on a
semigroupoid $G$; i.e., a partial order $\preceq$ such that for any
$a\in G$, the set $\{b\in G : b\preceq a\}$ is finite.  Since
$\preceq$ is a partial order, it is permissible to use the notation
$y\prec x$ to mean $y\preceq x$, but $y\ne x$.

As before, a set $F$ is lower if for all $a\in F$, $\{b: b\preceq a\}
\subset F$.  Clearly the intersection of lower sets is again lower.
For $z\in G$, let $S_z=\{x\preceq z\}$.  This is a lower set.
Furthermore, if $H$ is any subset of $G$, $x$ is minimal in $H$ is
equivalent to $S_x \cap H = \{x\}$.  By the artinian assumption, $S_x$
is finite.  Note that $b\preceq a$ is equivalent to $S_b \subseteq
S_a$.  (In fact there is an equivalence between artinian partial
orders $(G,\preceq)$ and functions $\lambda:G\to {\mathfrak S}_G$,
${\mathfrak S}_G$ the set of all finite subsets of $G$, $\lambda$
injective and $\lambda(\lambda(x)) = \lambda(x)$, where $a\preceq b$
if and only if $\lambda(a) \subseteq \lambda(b)$.)

\begin{lemma}
  \label{lem:minimal}
  Each nonempty subset $H$ of $G$ contains a minimal element with
  respect to an artinian partial order $\preceq$.
\end{lemma}

\begin{proof}
  Clearly any finite subset has a least element.  Suppose $H$ is any
  nonempty set and choose $z \in H$.  Now $S_z\cap H$ is a nonempty
  finite set, so it has a minimal element $x\in H$.  Since
  $S_x\subseteq S_z$, $S_x\cap H = S_x\cap (S_z\cap H) = \{x\}$; that
  is, $x$ is minimal in $H$.
\end{proof}

For a semigroupoid $G$ with artinian order $\preceq$, we define a
\textit{stratification} of $G$ as follows.  Set $G_0 = G_e$.  For
natural numbers $n$, define
\begin{equation*}
  G_{n} = \{ x\in G : y \prec x \Rightarrow y\in G_m\text{ for some }m
  < n \text{ and }y \prec x \text{ for some }y\in G_{n-1}\}.
\end{equation*}
We call $G_n$ the \textit{$n^{\mathrm{th}}$ stratum} with respect to
the order $\preceq$ and $\{G_n\}$ where $G_n$ is nonempty \textit{a
  stratification of $G$} with respect to the order $\preceq$.

\begin{lemma}
  \label{lem:stratified}
  For every $g\in G$ there is a unique $n\in\mathbb N$ such that $g\in
  G_n$.
\end{lemma}

\begin{proof}
  Suppose $H=G\setminus \cup_{0}^\infty G_n$ is nonempty.
  From Lemma \ref{lem:minimal}, $H$ has a minimal element $z$ (with
  respect to~$\preceq$).  In particular, $z\in S_z\subset
  \bigcup_{0}^\infty G_n$, a contradiction.
\end{proof}

For a lower set $F \subset G$ with respect to the order $\preceq$ we
define the \textit{stratification $\{F_n\}$ of $F$} with strata $F_n =
F \cap G_n$ where $F_n\neq \emptyset$.

The order $\leq$ which we originally introduced on semigroupoids
(where $b\leq a$ if and only if $a=zbw$ for some $z,w\in G$) is
artinian by definition of a semigroupoid.  Hence the above lemmas
apply to $G$ with this order.  There is another artinian order which
will be useful in proving the existence of inverses with respect to
the $\st$-product.

Define the {\it left order $\le_\ell$} on $G$ by declaring $y\leftless
x$ if there is an $a$ so that $x=ay$.

\begin{lemma}
  \label{lem:leftorder}
  The relation $\leftless$ is a partial order on $G$ which is more
  restrictive than the order $\leq$ on $G$; that is, if $y\leftless
  x$, then $y\le x$.
\end{lemma}

\begin{proof}
  The existence of idempotents implies that $x\leftless x$.  If
  $z\leftless y \leftless x$, then there exist $a,b$ so that $x=ay$
  and $y=bz$.  Hence, $x=a(bz)=(ab)z$ by the associative law and thus
  $z\leftless x$.  Finally, choosing $x=z$ above gives $x=(ab)x$.  By
  the strong idempotent law, it follows that $ab$ is idempotent; and
  then by nonexistence of inverses $a=b=e$ where $e$ is the idempotent
  so that $ex=x$.  Thus $x=ey$.  But by what it means to be
  idempotent, $ey=y$.  Hence if $x\leftless y$ and $y\leftless x$,
  then $x=y$.  This proves that $\leftless$ is an order on $G$.

  If $y\leftless x$, then $x=ay$ for some $a$.  There is always an
  idempotent $f$ so that $xf=x$.  Thus, $x=xf= ayf$ (by associativity)
  and $y\le x$.  Hence $\{y: y\leftless x\} \subseteq \{y: y\leq x\}$.
  The latter set is finite, so both are finite.
\end{proof}  

We use the notation $\{F_n^\ell\}$ for the stratification of a lower
set $F$ with respect to the left order, and for $z\in G$, we write
$S_z^\ell$ for $S_z$ with respect to the left order.

\subsection{$\st$-inverses}
\label{subsec:inverses}

We next prove the statement about inverses of matrices with respect to
the $\st$-product made in the introduction.  A similar (and in fact
easier) proof works for inverses of functions with respect to the
$\st$-product.  The arguments in the proof also apply to matrices over
any $C^*$-algebra, though the theorem is stated for matrices over
$\mathbb C$.

\begin{theorem}
  \label{thm:inverse}
  Let $F$ be a lower set.  A matrix $A\in M(F)$ is $\st$-invertible if and
  only if $A_{ab}$ is invertible for all $a,b\in F_e$.  Furthermore
  the inverse is unique.
\end{theorem}

\begin{proof}
  For $a,b\in F_e = F\cap G_e$, the term $(B\s A) (a,b)$ is
  essentially the Schur product (that is, $(B\s A) (a,b) =
  B_{ab}A_{ab}$ if $a,b\in F_e$) and if $B$ is a $\st$-inverse of $A$,
  then $(B\s A)(a,b)=1$ in this case.  Thus, $A_{ab}$ is invertible.

  The proof of the converse proceeds as follows.  Under the hypotheses
  of the theorem, a left $\st$-inverse $B$ for $A$ is constructed which
  itself satisfies the hypotheses of the theorem.  By what has already
  been proved, $B$ then has a left $\st$-inverse $C$.  Associativity of
  the $\st$-product guarantees that $C=A$ and thus $B$ is also a right
  $\st$-inverse for $A$.  Uniqueness of the $\st$-inverse similarly
  follows from the construction.
 
  So assume $A_{ab}$ is invertible for all $a,b\in F_0^\ell= F\cap
  G_e$.  Let $\{F_n^\ell\}$ be the left stratification of $F$.  Define
  $P_{jk} = \{(a,b) : a\in F_j^\ell,\ b\in F_k^\ell\}$, and $Q_N =
  \bigcup_{j,k \le N} P_{jk}$.  The proof proceeds by induction on
  $N$.

  We require
  \begin{equation}
    \label{eq:1}
    (B\s A)(a,b) = \sum_{pq=a}\sum_{rs=b} B_{pr}A_{qs} =
    \begin{cases}
      1 & \text{both }a,b\in G_e \\
      0 & \text{otherwise.}
    \end{cases}
  \end{equation}
  In the case $N=0$ (so that $a,b\in F\cap G_e$), the choice
  $B_{ab}=A_{ab}^{-1}$ is the unique solution to this equation.
  
  Now suppose that $B_{ab}$ have been defined for $(a,b)\in Q_N$
  satisfying equation (\ref{eq:1}) and suppose $(a,b)\in
  Q_{N+1}\backslash Q_N$.  Isolating the $(a,b)$ term in equation
  (\ref{eq:1}) gives
  \begin{equation*}
    0 = B_{ab}A_{ef} + \sum_{\scriptstyle pq=a\atop\scriptstyle p\neq
      a}\sum_{\scriptstyle rs=b\atop\scriptstyle r\neq b} B_{pr}A_{qs},
  \end{equation*}
  where $e,f\in G_e$ with $ae=a$ and $bf=b$.  In the second term on
  the right hand side, $p\leftsless a$ and $r\leftsless b$.  In
  particular, $p,r\in Q_N$ and the matrix $B_{pr}$ is already defined.
  Since $A_{ef}$ is invertible, $B_{ab}$ is uniquely determined.
\end{proof}

\begin{lemma}
  \label{lem:restricted_inverse}
  Let $L, F$ be lower sets in $G$ with $L\supset F$.  Suppose $A\in
  M(L)$ is $\st$-invertible.  Then $A|_F$ is $\st$-invertible and
  $(A|_F)^{-1_\st} = A^{-1_\st}|_F$.
\end{lemma}

\begin{proof}
  This follows by observing that $A|_F \s A^{-1_\st}|_F = (A\s
  A^{-1_\st})|_F = [1]_F$.
\end{proof}


\section{Reproducing kernels}
\label{sec:reproducing-kernels}

\subsection{Generalized Szeg\H{o} kernels}
\label{sec:gener-toepl-oper}

In this section we investigate those kernels which play the role
over semigroupoids of Szeg\H{o} kernels.  Recall, for a function
$\varphi$ defined on a lower set $F$, the $n$-fold $\st$-product of
$\varphi$ with itself is denoted $\varphi^{n_\st}$.  We use $A^{n_\st}$
similarly when $A$ is a matrix.

\begin{theorem}
  \label{thm:pos_inv2}
  Let $A\in M(F)$ be positive, and suppose $\|A^{n_\st}\| \to 0$ as
  $n\to\infty$.  Then $[1]-A$ is invertible $($with respect to the
  $\st$-product$)$ and $([1]-A)^{-1_\st} \geq 0$.  In particular, the
  result holds if $\|A\| < 1$.
\end{theorem}

\begin{proof} Observe that under the hypotheses, $[A^{n_\st}]_{e,e} =
  A^n_{e,e} \to 0$ as $n\to\infty$ for $e\in F_e$.  Hence $|A_{e,e}| <
  1$ for all $e\in F_e$.  Positivity of $A$ then implies that
  $|A_{e,f}| < 1$ for all $e,f\in F_e$.  Consequently $[1]-A$ is
  invertible.

  It is easily seen that
  \begin{equation*}
    1+A+A^{2_\st} + \cdots + A^{n_\st} = ([1]-A)^{-1_\st}\s ([1] -
    A^{(n+1)_\st}).
  \end{equation*}
  But $[1] - A^{(n+1)_\st} \to [1]$ as $n\to\infty$ and $1+A+A^{2_\st} +
  \cdots + A^{n_\st}$ is an increasing sequence of positive operators,
  and so converges strongly to $([1]-A)^{-1_\st}$.  Thus
  $([1]-A)^{-1_\st} \geq 0$.  The last part of the theorem follows from the
  submultiplicativity of the operator norm.
\end{proof}

It is not difficult to verify that the above arguments also work if
we instead consider matrices over a unital $C^*$-algebra.

\begin{corollary}
  \label{cor:pos-inv2}
  If $\mathfrak C$ is a unital $C^*$-algebra, $F$ a finite lower set,
  and $\varphi\in P(F,\mathfrak C)$.  Suppose that
  \begin{equation}
    \label{eq:extra}
    \lim_{n\to \infty} \varphi^{n}(a) =0
  \end{equation}
  for each $a\in F_e$.  Then $([1]-\varphi\varphi^*)^{-1_\st} \in
  M(F,\mathfrak C)$ is well defined and positive.  In particular, if
  $\|\sT(\varphi)\| < 1$, the result follows (and in this case $F$
  need not be finite).
\end{corollary}

\begin{proof}
  Let $A(a,b) = \varphi(a)\varphi(b)^*$.  Since for $e\in F_e$,
  $\|A^{n_\st}(e,e)\| = \|A^n(e,e)\| \to 0$ as $n\to\infty$, it follows
  that $\|A(e,f)\| < 1$ for all $e,f\in F_e$.  A counting argument in
  the same vein as that following Definition~\ref{def:test-functions}
  then shows that $\|A^{n_\st}(a,b)\| \to 0$ as $n\to\infty$, and so
  since $F$ is assumed to be finite, $\lim_{n\to\infty}\|A^{n_\st}\| =
  0$.  The conditions for Theorem~\ref{thm:pos_inv2} hold and the
  result follows directly.

  If $\|\sT(\varphi)\| < 1$ then since
  \begin{eqnarray*}
    (A\s 1)_{a,b} &=& \sum_{pq=a}\sum_{rs=b}
    \varphi(p)\varphi(r)^*\ip{1_q}{1_s}\\
    &=& \sum_{pq=a}\sum_{rq=b}\varphi(p)\varphi(r)^*
  \end{eqnarray*}
  and
  \begin{eqnarray*}
    (\sT(\varphi)\sT(\varphi)^*)_{a,b} &=& \sum_q
    [\sT(\varphi)]_{a,q}[\sT(\varphi)^*]_{q,b}\\
    &=& \sum_q [\sT(\varphi)]_{a,q}[\sT(\varphi)]_{b,q}^*\\
    &=& \sum_{pq=a}\sum_{rq=b}\varphi(p)\varphi(r)^*.
  \end{eqnarray*}
  the result is then a consequence of the last statement of
  Theorem~\ref{thm:pos_inv2}, since $\|A\s 1\| = \|A\|_{cb}$ dominates
  the operator norm.
\end{proof}

In what follows the theorem will be applied to test functions $\psi$
and more generally the evaluation $E$.  That $E$ satisfies the
hypothesis of Corollary~\ref{cor:pos-inv2} is equivalent to item (ii)
in the definition of test functions
(Subsection~\ref{subsec:intro-test-functions}).

\begin{lemma}
 \label{lem:pos-inv-E}
 For each $a\in G_e$, the sequence $E^{n_\st}(a)$ from $\CT$ converges
 to $0$.
\end{lemma}

\subsection{The multiplier algebra for a single kernel}
\label{sec:mult-algebra-single}

Let $\mathbf{k}:G\times G\to \mathbb C$ be a positive kernel.  For $b\in G$,
the function $k_b:G\to\mathbb C$ defined by $k_b = \mathbf{k}(\cdot,b)$ is
point evaluation at $b$.  In the usual way we form a sesquilinear form
$\ip{\cdot}{\cdot}$ on linear combinations of kernel functions by
setting $\ip{k_b}{k_a} = \mathbf{k}(a,b)$ and modding out by the kernel.  We
then complete to get a Hilbert space, $H^2(\mathbf{k})$.

On $H^2(\mathbf{k})$ addition is defined term-wise.  The multiplier
algebra $H^\infty(\mathbf{k})$ consists of the collection of operators
$T_\varphi: f\mapsto \varphi \star f$ for functions $\varphi:G\to
\mathbb C$ satisfying $\varphi \s f\in H^2(\mathbf{k})$ for each $f\in
H^2(\mathbf{k})$.  (The product is well defined by the assumption that
$G$ is artinian.)  Note that $H^\infty(\mathbf{k})$ is nonempty, since
it contains $T_{\delta}$, $\delta$ the $\st$-product identity for
functions on $G$.  The closed graph theorem implies that the elements
of $H^\infty(\mathbf{k})$ are bounded.

Observe that for $f\in H^2(\mathbf{k})$,
\begin{eqnarray*}
  \ip{T_\varphi f}{k_a} &=& (\varphi\s f)(a) \\
  &=& \sum_{bc=a} \varphi(b)f(c) \\
  &=& \sum_{bc=a} \varphi(b) \ip{f}{k_c} \\
  &=&  \sum_{bc=a} \ip{f}{\varphi(b)^* k_c} \\
  &=& \ip{f}{\sum_{bc=a} \varphi(b)^* k_c},
\end{eqnarray*}
which gives the formula $T_\varphi^* k_a = \sum_{bc=a} \varphi(b)^*
k_c$.

For a lower set $F$, if we set ${\mathcal M}(F)$ to the closed linear
span of kernel functions $k_a$, $a\in F$, then the usual sort of
argument gives ${\mathcal M}(F)$ invariant for the adjoints of
multipliers $T_\varphi$.

The $\st$-product is useful in characterizing multipliers.  Indeed,
$\|T_\varphi^*|\mathcal M(F)\| \leq C$ is equivalent to $0 \leq
(\ip{(C^2-T_\varphi T_\varphi^*) k_a}{k_b})$ which by the previous
calculation is
\begin{equation*}
  \| T_\varphi^*|\mathcal M(F)\| \leq C \qquad\Longleftrightarrow
  \qquad C^2\mathbf{k}-\varphi^*\s \mathbf{k} \sh \varphi \geq 0.
\end{equation*}
In the above $\varphi^* \s \mathbf{k} \sh \varphi$ stands for
$(\varphi^*\s k) (k^*\sh \varphi)$ where $k(x) = k_x$ in the
factorization $\mathbf{k}(x,y) = k_xk_y^*$ for $x,y \in F$.

\subsection{The Toeplitz kernel}

A special case of interest is the kernel $s:F\times F\to \mathbb C$
given by $s(x,y)=1$ if $x=y$ and $0$ if $x\ne y$.  This kernel is
evidently positive and, as noted earlier is referred to as the
Toeplitz kernel.  It arises naturally by declaring $\ip{x}{y} =s(x,y)$
for $x,y\in F$ and extending by linearity.  That is, the Hilbert space
$H^2(s)$ is nothing more than the Hilbert space with orthonormal basis
indexed by $F$; i.e., $\mathbb C^F$.  The Toeplitz representation of
$\varphi:F\to \mathbb C$ determined by $s$ as in the previous
subsection is thus the Toeplitz representation $\sT (\varphi)$ of
$\varphi$.

Note that $\begin{pmatrix}s(x,y)\end{pmatrix}_{x,y\in F} =1\in M(F)$,
the usual identity matrix.

\subsection{Kernels and representations}
\label{subsec:ker-rep}

The results of Subsection \ref{sec:mult-algebra-single} have an
alternate interpretation.  Let $\mathbf{k}$ be a reproducing kernel on
$G$.  Recall that we use $P(F)$ to denote the complex valued functions
on the finite lower set $F\subset G$, which under the $\st$-product is
an algebra.  If we now compress $\mathbf{k}$ to $F$, it is still a
positive kernel (on $F$) which we continue to call $\mathbf{k}$ (or
$\mathbf{k}_F$ if it is not absolutely clear from the context).
Furthermore, since $F$ is lower
\begin{equation}
 \label{eq:mathK-F-Psi}
  F\times F \ni (x,y)\mapsto (([1]-\psi \psi^*)\s \mathbf{k})(x,y)
\end{equation}
is positive for each $\psi \in\PrF$.  In this case, any $\varphi\in
P(F)$ is a multiplier of $H^2(\mathbf{k}_F)$ since the algebra
generated by $\PrF$ is all of $P(F)$ and $\pi: P(F) \to
B(H^2(\mathbf{k}_F))$ defined by $\pi(\varphi) = \sT^F(\varphi)$ is a
representation of $P(F)$.  Further, the assumption that
(\ref{eq:mathK-F-Psi}) is positive implies $\|\pi(\psi)\|\le 1$ for
each $\psi \in \PrF$.
 
Define the functions
\begin{equation*}
  \chi_a (x) =
  \begin{cases}
    1 & x=a,\\
    0 & \text{otherwise}.
  \end{cases}
\end{equation*}
Routine calculation verifies $(\chi_a\s \chi_b)(x) = \sum_{pq=x}
\chi_a(p)\chi_b(q) = \chi_{ab}(x)$, where for convenience we take
$\chi_{ab} = 0$ if the product $ab$ is not in our partial
multiplication.

Clearly the set $\{\pi(\chi_a)\delta\}$ forms a spanning set for
$H^2(\mathbf{k}),$ and, since $\pi(\varphi) \delta =\varphi= 0$ if and
only if $\varphi =0$, it is in fact a basis.  Indeed it is a dual
basis to $\{k_a\}$, since
\begin{equation*}
  \ip{k_a}{\pi(\chi_b)\delta} = \ip{\sT^F(\chi_b)^* k_a}{\delta} =
  \sum_{pq=a} \ip{\chi_b^* (p) k_q}{\delta} = 
  \begin{cases}
    1 & \text{if }p=b=a,\ q\in F_e \\
    0 & \text{otherwise.}
  \end{cases}
\end{equation*}

In some cases it is possible to reverse the above, obtaining a kernel
from a representation $\mu: P(F)\to \mathcal B(\mathcal H)$.  For
instance, suppose $F$ is a finite lower set and assume that $\mu$ is
cyclic with dimension equal to the cardinality of $F$.  Write $\gamma$
for the cyclic vector for $\mu$, so that $\mathcal H$ is spanned by
$\{\ell_a = \mu(\chi_a)\gamma : a\in F\}$.  Since by assumption the
dimension of $\mu$ is the cardinality of $F$, this set is in fact a
basis for $F$.

If $\mu$ is to come from a kernel $\mathbf{k}$, we require that for any
function $\varphi$ on $F$,
\begin{equation*}
  \mu(\varphi)^* k_a = \sum_{pq=a} \varphi(p)^* k_q.
\end{equation*}
It suffices to have this for the functions $\chi_b$, in which case we
need
\begin{equation*}
  \mu(\chi_b)^* k_a = \sum_{pq=a} \chi_b(p) k_q = \sum_{bq=a} k_q.
\end{equation*}

Choose $\{k_a: a\in F\}$ to be a dual basis to $\ell_a$.  Then
compute,
\begin{equation*}
 \begin{split}
  \ip{\mu(\chi_b)^* k_a}{\ell_c}
    =\,& \ip{k_a}{\mu(\chi_b)\ell_c} \\
    =\,& \ip{k_a}{\mu(\chi_b)\mu(\chi_c)\gamma} \\
    =\,& \ip{k_a}{\mu(\chi_b \s \chi_c)\gamma} \\
    =\,& \ip{k_a}{\ell_{bc}} \\
    =\,& \begin{cases} 1 & bc=a \\ 0 & bc \ne a \end{cases} \\
    =\,& \ip{\sum_{bq=a} k_q}{\ell_c} .
 \end{split}
\end{equation*}
Since this is true for all $c\in F$, it follows that
\begin{equation*}
  \mu(\chi_b)^*k_a=\sum_{bq=c} k_q,
\end{equation*}
as desired.

It is worth considering the example where $\mu(\varphi) =
\sT(\varphi)$, the Toeplitz representation.  The function $\delta(x)$,
which is $1$ if $x\in F_e$ and zero otherwise is a cyclic vector for
$\mu$.  Moreover,
\begin{equation*}
  \mu(\chi_a)(x,y) = \sT(\chi_a)(x,y) = \sum_{py=x} \chi_{a}(p) =
  \chi_{ay}(x).
\end{equation*}
 Therefore
\begin{equation*}
  \ell_a(x) = \sum_{y} \mu(\chi_a)(x,y) \delta(y) = \sum_{y}
  \chi_{ay}(x) \delta(y) = \chi_a(x),
\end{equation*}
which is just the standard basis, and so the assumption that
$\{\ell_a: a\in F\}$ is a basis is automatically met.  In this case we
choose $k_a = \chi_a$, and the kernel is the Toeplitz kernel $s$.

If $F$ is infinite, this construction fails, since it need not be the
case that $\chi_a\in \HP$.

\subsection{$\mathbf{P(F)}$ as a normed algebra}

Given a finite lower set $F$, let $\pi_F:\HP \to P(F)$ denote the
mapping $\pi_F(\varphi)=\varphi|_{F}$.  The hypothesis on the
collection of test functions $\Psi$ imply that this mapping is onto
and so $\mbox{ker}(\pi_F)=\{\varphi\in \HP: \varphi|_F=0\}$.  Thus
$P(F)$ is naturally identified with the quotient of $\HP$ by
$\mbox{ker}(\pi_F)$ and this gives $P(F)$ a norm for which $\pi_F$ is
contractive.  There is an alternate candidate for a norm on $P(F)$
constructed in much the same way as the norm on $\HP$, called the
$\HPF$-norm.  Let $\mathcal K_\Psi^F$ denote the kernels $\mathbf{k}$
defined on $F$ for which
\begin{equation*}
  F\times F \ni (x,y)\mapsto  (([1]-\psi|_F \psi|_F^*)\s \mathbf{k})(x,y)
\end{equation*}
is a positive kernel and, for $\varphi \in P(F)$, say that
$\|\varphi\|\le C$ (here $C\geq 0$) provided for each
$\mathbf{k}\in\mathcal K_\Psi^F$, the kernel
\begin{equation*}
  F\times F \ni (x,y)\mapsto ((C^2[1]-\varphi \varphi^*)\s \mathbf{k})(x,y)
\end{equation*}
is positive.

The following lemma ultimately implies that the quotient norm
dominates the $\HPF$-norm.  Theorem \ref{thm:main-interpolate} then
says that these norms are the same.

\begin{lemma}
  \label{lem:reptok}
  Suppose $\mu:P(F)\to B(\mathcal H)$ is a cyclic unital
  representation of the finite lower set $F$ and let $\varphi\in \HP$
  be given.  Let $\pi_F : \HP \to P(F)$ be the restriction map, $\mu_F
  = \mu\circ\pi_F$.  If $\|\mu_F(\psi)\|\le 1$ for each $\psi \in
  \Psi$, but $\|\mu_F(\varphi)\|>1$, then there exists a
  $\mathbf{k}\in\mathcal K_\Psi$ so that the kernel
  \begin{equation*}
    F\times F \ni (x,y)\mapsto (([1]-\varphi\varphi^*)\s \mathbf{k})(x,y)
  \end{equation*}
  is not positive.  In particular, $\|\varphi\|>1$.
\end{lemma}

\begin{proof}
  Let $\gamma$ denote a cyclic vector for the representation $\mu$.
  Choose $f\in P(F)$ so that $\|\mu(f)\gamma\|=1$ but
  $\|\mu_F(\varphi)\mu(f)\gamma\|=1+\eta >1$, and $\epsilon$ so that
  $(1+\epsilon^2 \|f\|^2_{H^2(s)})(1+\eta/2)^2 = (1+\eta)^2$, where
  $\|f\|_{H^2(s)}$ is the norm of $f$ in the space with the Toeplitz
  kernel $s$.
  
  Recall that for a finite lower set $L$, $\sT^L$ denotes the Toeplitz
  representation with its cyclic vector $\delta^L$.  If $L
  \supseteq F$, let $\pi^L_F$ be the restriction of $P(L)$ to $P(F)$
  and set $\mu^L_F = \mu\circ\pi^L_F$.  As above, define $\pi_L : \HP
  \to P(L)$ to be the restriction map.
  
  For $L \supseteq F$ lower, there is a finite dimensional Hilbert
  space given by $\mathcal H_L = \{ \mu^L_F(h)\gamma \oplus \epsilon
  \sT^L(h)\delta^L: h\in P(L)\}$.  Recall $\sT^L(h)\delta^L=h$.
  Define a representation $\rho_L : P(L) \to B(\mathcal H_L)$ by
  \begin{equation*}
    \rho_L(g)(\mu^L_F(h)\gamma \oplus \epsilon\sT^L(h)\delta^L) = 
    \mu^L_F(g\s h)\gamma \oplus \epsilon\, g\s h
  \end{equation*}
  Since for $\psi \in \Psi$, $h\in P(L)$,
  \begin{equation*}
    \begin{split}
      \|\rho_L(\pi_L(\psi)) (\mu(h)\gamma \oplus \epsilon  h)\|
      &= \|\mu(\pi_F (\psi)) \mu(h)\gamma \oplus 
           \sT^L(\pi_L\psi) \epsilon h\| \\
      &\le \max\{\|\mu(\pi_F(\psi))\|, \|\sT^L(\pi_L\psi)\|\}
      \|\mu(h)\gamma \oplus \epsilon  h\| \\
      &\le \|\mu(h)\gamma \oplus \epsilon  h\|,
    \end{split}
  \end{equation*}
  $\|\rho_L(\psi)\|\le 1$, and in particular, taking $L=F$ we
  have $\|\rho_F(\psi)\|\le 1$.

  From the discussion in Subsection \ref{subsec:ker-rep}, there is a
  kernel $\mathbf{k}^F$ on $F$ which implements the representation
  $\rho_F$.  In particular, since
  \begin{equation*}
    \begin{split}
      {\left\| \rho_F (\varphi)\frac{\mu(f)\gamma \oplus \epsilon
            f}{\sqrt{1+\epsilon^2 \|f\|^2_{H^2(s)}}}\right\|}^2 & =
      \frac{1}{1+\epsilon^2 \|f\|^2_{H^2(s)}}
      \left(\|\mu(\varphi)\mu(f)\gamma\|^2 +
        \epsilon^2 \| \varphi \s f \|^2\right) \\
      & \geq \frac{1}{1+\epsilon^2 \|f\|^2_{H^2(s)}} (1+\eta)^2
      \\
      & = (1+\eta/2)^2,
    \end{split}
  \end{equation*}
  the kernel 
  \begin{equation*}
    F\times F\ni  (x,y) \mapsto (([1]-\varphi\varphi^*)\s \mathbf{k}^F)(x,y)
  \end{equation*} 
  is not positive.
  
  Define $\mathbf{k}:G\times G\to \mathbb C$ by 
  \begin{equation*}
    \mathbf{k}(a,b)=\begin{cases} \mathbf{k}^F(a,b) & \mbox{ if }
      (a,b)\in F\times F \\ 
      \frac{1}{\epsilon^2} s(a,b) & \mbox{ if } (a,b)\notin F\times F.
    \end{cases}
  \end{equation*}
  In particular, if $a\in F$ and $b\notin F$ (or vice-versa), then
  $\mathbf{k}(a,b)=0$.  We will complete the proof by showing
  $\mathbf{k}\in\mathcal K_\Psi$.
  
  The representation $\rho_L$ defined as above is non-degenerate for
  any $L \supseteq F$, in the sense of the discussion in Subsection
  \ref{subsec:ker-rep}.  In particular, for any such $L$ there is a
  reproducing kernel $\mathbf{k}^L$ which implements this
  representation.  Consequently, for each $\psi\in \Psi$,
  \begin{equation*}
    L\times L \ni (x,y)\mapsto (([1]-\psi\psi^*)\s \mathbf{k}^L)(x,y)
  \end{equation*}
  is positive.  Our goal now is to show that
  $\mathbf{k}^L(x,y)=\mathbf{k}(x,y)$ for $x,y\in L$ from which it
  will follow that $\mathbf{k}\in\mathcal K_\Psi$.

  For this, once again recall the construction of $\mathbf{k}^L$ from
  $\rho_L$.  Let $\ell_a^L=\rho_L(\chi_a)h$, where $h=\gamma\oplus
  \epsilon \delta^L= \mu^L_F(\delta^L)\gamma \oplus \epsilon
  \sT^L(\delta^L)\delta^L$ is the cyclic vector for the representation
  $\rho_L$.  Next, let $k_b^L$ denote a dual basis to the basis
  $\ell_a^L$ and define $\mathbf{k}^L$ by $\mathbf{k}^L(a,b) =
  \ip{k^L_b}{k^L_a}$.
 
  We calculate
  \begin{equation*}
    \ell_a^L = \rho_L(\chi_a)h = \mu^L_F(\chi_a\s \delta^L)\gamma
    \oplus \epsilon\, \chi_a\s \delta^L = \mu^L_F(\chi_a) \oplus
    \epsilon\, \chi_a,
  \end{equation*}
  which reduces to $\{0\} \oplus \epsilon\, \chi_a$ if $a\notin F$,
  and which equals $\ell_a^F \oplus \{0\}$ if $a\in F$.  Hence the
  dual basis is
  \begin{equation*}
    k_a^L =
    \begin{cases}
      k_a^F \oplus \{0\} & a\in F,\\
      \{0\} \oplus \epsilon\, \chi_a & \text{otherwise,}
    \end{cases}
  \end{equation*}
  via which we immediately verify that $\mathbf{k}^L = \mathbf{k}$.
\end{proof}

\subsection{ Toeplitz representation for $C^*$-algebra-valued
  functions}
\label{subsec:more-standard-Toeplitz}

The notion of the Toeplitz representation naturally generalizes to
functions $f:F\to \mathfrak C$, where $F$ is a lower set and
$\mathfrak C$ is a $C^*$-algebra with $[\sT (\varphi)]_{a,b}\in
\mathfrak C$ and $\sT(\varphi) \in M(F,\mathfrak C)$, the $\mathfrak
C$-valued matrices labelled by elements of $F$.

\begin{lemma}
  \label{lem:reptoeplitzrep}  
  Suppose that $\mathcal C$ is another $C^*$-algebra.  If
  $\rho:\mathfrak C\to \mathcal C$ is a unital $*$-representation,
  then
  \begin{equation*}
    (1\otimes \rho) (\sT(f))=\sT(\rho \circ f)
  \end{equation*}
  and moreover,
  $\|\sT(f)\|\ge \|\sT(\rho\circ f)\|.$
\end{lemma}

\begin{proof}
  Simply compute
  \begin{equation*}
    \begin{split}
      (1\otimes \rho) (\sT(f))= &[ \rho([\sT(f)]_{a,b})]_{a,b} \\
      =\,& \Bigg[\rho \Bigg(\sum_{cb=a} f(c)\Bigg)\Bigg]_{a,b}\\
      =\,&\Bigg[\sum_{cb=a} \rho(f(c))\Bigg]_{a,b}\\
      =\,& [[\sT(\rho\circ f)]_{a,b}]_{a,b}.
    \end{split}
  \end{equation*}
  The norm estimate follows since $\rho$ is completely contractive.
\end{proof} 

In our applications of this lemma $f$ will be the function $E:F\to
\CT$ and $\rho:\CT \to B(\mathcal E)$ will be the representation
arising in a $\CT$-unitary colligation.


\section{Factorization}
\label{sec:factorization}

\begin{proposition}
 \label{prop:factorization}
 If $\Gamma:G\times G \to \CT^*$ is positive, then there exists a
 Hilbert space $\mathcal E$ and a function $L:G\to B(\CT, \mathcal E)$
 such that
 \begin{equation*}
   \Gamma(x,y)(fg^*) = \ip{L(x)f}{L(y)g}
 \end{equation*}
 for all $f,g\in \CT$.
  
  Further, there exists a unital $*$-representation $\rho:\CT\to
  B(\mathcal E)$ such that $L(x)ab=\rho(a) L(x)b$ for all $x\in G$,
  $a,b\in\CT$.
\end{proposition}

\begin{proof}
  The proof is a variant on a usual proof of the factorization of
  positive semidefinite kernels. See the book \cite{MR1882259}
  Theorem~2.53, Proof~1. The statement  should be compared with
  a similar result in \cite{MR2106336}.  
  
  Let $W$ denote a vector space with basis labelled by $G$.  On the
  vector space $W\otimes \CT$ introduce the positive semidefinite
  sesquilinear form induced from
  \begin{equation*}
    \ip{x \otimes f}{y\otimes g} =\Gamma(x,y)(fg^*),
  \end{equation*}
  where $x,y\in G$ and $f,g\in\CT$, making $W\otimes \CT$ into a
  pre-Hilbert space which is made into the Hilbert space $\mathcal E$
  by the standard modding out and completion.
  
  One verifies that this is indeed positive as a consequence of the
  hypothesis that $\Gamma$ is positive.  Define $L(x)a= x\otimes a$.
  Since for $a\in\CT$,
  \begin{equation*}
    \begin{split}
      \|L(x)a\|^2 = & \ip{L(x)a}{L(x)a} \\
      =\,& \Gamma(x,x)(a^* a) \\
      \le \,&  \|\Gamma(x,x)\| \|a^* a\| 
    \end{split}
  \end{equation*}
  $L(x)$ does indeed define a bounded operator on $\CT$ with
  $\|L(x)\|^2\le \|\Gamma(x,x)\|$.
  
  As for the $*$-representation, it is induced by the left regular
  representation of $\CT$.  That is, define $\rho:\CT \to B(\mathcal
  E)$ by $\rho(a)(x\otimes f)= x\otimes af$.  To see that this is
  indeed bounded, first note that $\|a\|^2-a^*a$ is positive
  semidefinite in $\CT$ and hence there exists a $b$ so that
  $\|a\|^2-a^*a = b^*b$.  Thus,
  \begin{equation*}
    \begin{split}
      \|a\|^2 &\Big\|\sum x_j\otimes f_j\Big\|^2 - \Big \| \sum x_j
      \otimes af_j \Big\|^2  \\
      =\,& \|a\|^2 \sum \Gamma(x_j,x_\ell)(f_\ell^* f_j)-\sum
      \Gamma(x_j,x_\ell)(f_\ell^* a^* a  f_j) \\
      =\,& \sum \Gamma(x_j,x_\ell)(f_\ell^* b^* bf_j) \ge 0
    \end{split}
  \end{equation*}
  where the inequality is a result of the assumption that $\Gamma$ is
  positive.  This shows at the same time that $\rho$ is well defined.

  We also have that $\rho$ is unital, since $\rho(1)(x\otimes f)=
  x\otimes 1f=x\otimes f$.
  
  Finally,
  \begin{equation*}
    \begin{split}
      \ip{\rho(a^*)(x \otimes f)}{ y\otimes g}
      =\,& \ip{x \otimes a^*f}{y\otimes g} \\
      =\,&\Gamma(x,y)(g^* a^* f)\\
      =\,&\ip{x \otimes f}{y\otimes ag} \\
      =\,& \ip{x \otimes f}{ \rho(a) (y\otimes g)}\\
      =\,&\ip{\rho(a)^* (x \otimes f)}{ y\otimes g}\\
    \end{split}
  \end{equation*}
  so that $\rho(a^*)=\rho(a)^*$.
\end{proof}


\section{The cone $\mathcal C_F$ and compact convex set $\Phi_F$}

Given a finite subset $F\subset G$, let $M(F,\CT^*)^+$ denote the
collection of positive kernels $\Gamma:F\times F\to \CT^*$ and define
the cone
\begin{equation*}
  \mathcal C_F=\{ \begin{pmatrix} (\Gamma\sh ([1]-EE^*))(x,y)
  \end{pmatrix}_{x,y\in F}: \Gamma \in M(F,\CT^*)^+\}.
\end{equation*}

\subsection{The cone is closed}
\label{subsec:closedcone}

\begin{theorem}
  \label{thm:closedcone}
  Let $F$ be a finite lower set.  The cone $\mathcal C_F$ is closed in
  $M(F)$.
\end{theorem}

\begin{proof} Let $M = \Gamma\sh ([1]-EE^*) \in \mathcal C_F$, where
  $\Gamma:F\times F\to \CT^*$ is positive.  Positivity of $\Gamma$
  means in particular that if $\tilde F$ is a subset of $F$, and
  $\{f_q\}_{q\in \tilde F}$ is any collection of elements of $\CT$,
  then
  \begin{equation}
    \label{eq:7}
    \sum_{p,q\in \tilde F} \Gamma(p,q)f_pf_q^* \geq 0.
  \end{equation}
  
  For convenience we define $b \leq_r a$ to mean $bc = a$ for some
  $c$.  As in Lemma~\ref{lem:leftorder}, this can be shown to be an
  order on $G$ and $b \leq_r a$ implies $b\leq a$.

  Fix $x\in F$, and suppose $e$ is idempotent with $xe=x$.  Taking
  $F_x = \{y : y \leq_r x\}$, we get a finite subset of $F$.  For
  $q\in F_x$, set
  \begin{equation*}
    \tE(q) = \sum_{\scriptstyle p\atop\scriptstyle qp=x} E(p).
  \end{equation*}
  Observe that $\tE(x) = E(e)$.  With this notation, we have
  \begin{equation}
    \label{eq:8}
    \begin{split}
      M_{x,x} &\, = \Gamma(x,x)(1-\tE(x)\tE(x)^*) -
      \sum_{\scriptstyle q<_r x\atop\scriptstyle s<_r x}
      \Gamma(q,s)(\tE(q)\tE(s)^*) \\
      \qquad &- \sum_{q<_r x} \Gamma(q,x)(\tE(q)\tE(x)^*)
      - \sum_{s<_r x} \Gamma(x,s)(\tE(x)\tE(s)^*).
    \end{split}
  \end{equation}
  
  Choose $\epsilon > 0$ small enough that
  $1-(1+\epsilon^2)\tE(x)\tE(x)^* > 0$.  This can be done since
  $1-\tE(x)\tE(x)^* = 1-E(e)E(e)^* > 0$ by property (ii) of
  Definition~\ref{def:test-functions}.  Let $f_x = -\epsilon \tE(x)$,
  $f_q = (e^{-i\theta_q}/\epsilon) \tE(q)$ for $q\in F_x$, $q\neq
  x$, and $\theta_q = \arg (\sum_{q<_r x}\Gamma(q,x)\tE(q)\tE(x)^*)$.
  With this choice, \eqref{eq:7} gives
  \begin{equation}
    \label{eq:9}
    2|\sum_{q<_r x} \Gamma(q,x)(\tE(q)\tE(x)^*)| \leq (1/\epsilon)^2
    \sum_{q<_r x} \sum_{s<_r x} \Gamma(q,s)(\tE(q)\tE(s)^*) +
    \epsilon^2 \Gamma(x,x)(\tE(x)\tE(x)^*).
  \end{equation}
  Combining the inequality in \eqref{eq:9} with \eqref{eq:8} we have
  \begin{equation}
    \label{eq:10}
    \Gamma(x,x)(1-(1+\epsilon^2)\tE(x)\tE(x)^*) \leq M_{x,x} +
    \left(1+\tfrac{1}{\epsilon^2}\right) \sum_{q<_r x}
    \sum_{s<_r x} \Gamma(q,s)(\tE(q)\tE(s)^*).
  \end{equation}
  Furthermore, positivity of $\Gamma$ and a calculation as for
  \eqref{eq:9} yields for $g\in \CT$
  \begin{equation*}
    2|\Gamma(x,y) g| \leq \Gamma(x,x)\,1 + \Gamma(y,y)\,gg^* \leq
    \|\Gamma(x,x)\| + \|\Gamma(y,y)\|\,\|g\|^2,
  \end{equation*}
  and so
  \begin{equation}
    \label{eq:12}
    \|\Gamma(x,y)\| \leq \tfrac{1}{2}\left(\|\Gamma(x,x)\| +
    \|\Gamma(y,y)\|\right).
  \end{equation}
  
  We show by induction on (right) strata that for each $p,q\in F$,
  there is a constant $c_{p,q}$, independent of $\Gamma$, such that
  $\|\Gamma(p,q)\| \leq c_{p,q} \|M\|$.  By \eqref{eq:12}, it suffices
  to prove this for $p=q$.  Since $F$ is assumed finite, it will then
  follow that $\|\Gamma\| \leq c \|M\|$ for some $c\geq 0$ and
  independent of $\Gamma$.

  To begin with, if $e\in F$ is idempotent, then $M_{e,e} =
  \Gamma(e,e) (1-E(e)E(e)^*)$, and since $1-E(e)E(e)^* > 0$, we have
  that $c_{e,e}$ exists.  Now suppose that we have $c_{p,q}$ for all
  $p,q$ in the $(n-1)$st and lower strata.  Let $x$ be in the $n$th
  stratum.  Then by the induction hypothesis and \eqref{eq:10}, we
  find $c_{x,x}$.

  Let $\{M_j\}$ be a bounded sequence from $\mathcal C_F$,
  $M_j=\Gamma_j\sh ([1]-EE^*)$, so that
  \begin{equation*}
    M_j(x,y) = \Gamma_j(x,y)(1) -
    \sum_{pq=x}\sum_{rs=y}\Gamma_j(q,s)(E(p)E(r)^*),\quad x,y \in F.
  \end{equation*}
  Then $\{\Gamma_j\}$ is a bounded sequence in $M(F,\CT^*)^+$; i.e.,
  there is a uniform bound on the norm of the linear functional
  $\Gamma_j(x,y)$ independent of $x,y,j$.  It follows from weak-$*$
  compactness, that there exists $\Gamma \in M(F,\CT^*)$ and a
  subsequence $\{\Gamma_{j_\ell}\}$ of $\{\Gamma_j\}$ so that for each
  $x,y\in F$, the sequence $\{\Gamma_{j_\ell}(x,y)\}$ converges to
  $\Gamma(x,y)$ weak-$*$.  In particular,
  $\{\Gamma_{j_\ell}(p,r)(E(q)E(s)^*)\}$ converges to
  $\Gamma(p,r)(E(q)E(s)^*)$ for each $p,q,r,s$ (and also with
  $E(q)E(s)^*$ replaced by $1$).  If now $\{M_j\}$ converges to some
  $M$, then
  \begin{equation*}
    M=\lim_{\ell\to \infty}  \begin{pmatrix} \Gamma_{j_\ell}\sh
      ([1]-EE^*)(x,y) \end{pmatrix}_{x,y\in F}
    = \begin{pmatrix} \Gamma \sh ([1]-EE^*)(x,y) \end{pmatrix}_{x,y\in
      F}.
  \end{equation*}
   
  If $f:F\to \CT$, then
  \begin{equation*}
    0\le \sum_{x,y\in F} \Gamma_{j_\ell}(x,y)(f(x)f(y)^*) \rightarrow
    \sum_{x,y\in F} \Gamma(x,y)(f(x)f(y)^*),
  \end{equation*}
  which shows that $\Gamma$ is positive and completes the proof.
  \end{proof}

\subsection{The cone is big}

\begin{lemma}
  \label{lem:evalmu}
  Let $\TT$ be a set of test functions for $G$. 
   For each $\psi\in\TT$ the function $\Gamma_\psi:G\times G\to
  \CT^*$ given by
  \begin{equation*}
    \Gamma_\psi(x,y)(f)= ([1]-\psi\psi^*)^{-1_\st}(x,y)f(\psi), \qquad
    f\in \CT=C(\TT),
  \end{equation*}
  is a positive kernel.
\end{lemma}

\begin{proof}
  For each $x,y\in G$, the functional $\Gamma_{\psi}(x,y)$ is a
  multiple of evaluation at $\psi$ and hence does indeed define an
  element of $\CT^*$.

  For a finite lower set $F\subset G$ and a function $f:F\to \CT$,
  \begin{equation*}
   \begin{split}
     \sum_{x,y\in F} \Gamma_\psi(x,y)(f(x)f(y)^*)
     = \,& \sum_{x,y\in F} ([1]-\psi\psi^*)^{-1_\st}(x,y)(f(x)(\psi)
     f(y)(\psi)^*) \\
     =\,& \sum_{x,y\in F} ([1]-\psi\psi^*)^{-1_\st}(x,y)(g(x)g(y)^*)\\
   \end{split}
 \end{equation*}
 where $g:F\to \mathbb C$ is given by $g(x)=f(x)(\psi)$ and $g$ is the
 vector with $x$ entry $g(x)$.  By Corollary~\ref{cor:pos-inv2},
 \begin{equation*}
   F\times F \ni (x,y)\mapsto ([1]-\psi\psi^*)^{-1_\st}(x,y)
 \end{equation*}
 is a positive matrix in $M(F)$.  The conclusion follows.
\end{proof}

\begin{lemma}
  \label{lem:Kpos}
  Suppose $F\subset G$ is a finite lower set.  The cone $\mathcal C_F$
  contains all positive matrices.  In particular, it contains $[1]$
  and so has non-trivial interior.
\end{lemma}

\begin{proof}
  Let $\Gamma_\psi$ denote the positive kernel from the previous
  lemma.  Then
  \begin{equation*}
    [1](x,y)=\Gamma_\psi\sh ([1]-EE^*)(x,y) \in\mathcal C_F.
  \end{equation*}
  On the other hand, if $P\in M(F)$ with $P\geq 0$, then $P\sh
  \Gamma_\psi \geq 0$ and $P = P\sh [1] = P\sh(\Gamma_\psi\sh
  ([1]-EE^*))$.  Thus $\mathcal C_F$ contains all positive $P\in
  M(F)$.
\end{proof}

\begin{lemma} 
  \label{lem:conjcone}
  The cone $\mathcal C_F$ is closed under conjugation; i.e., if
  $M=(M(x,y))\in\mathcal C_F$ and $c:F\to \mathbb C$, then $c\s M \sh
  c^* \in\mathcal C_F$, where $(c\s M \sh c^*)(x,y) =
  \sum_{pq=x}\sum_{rs=y} c(p)M(q,s)c^*(r)$.
\end{lemma}

\begin{proof}
  If $M = \Gamma \sh ([1]-EE^*) \in \mathcal C_F$, then $c\s M \sh c^*
  = \tilde \Gamma \sh ([1]-EE^*)$, where $\tilde \Gamma = c\s \Gamma
  \sh c^* \geq 0$.
\end{proof}

\subsection{Separation}

\begin{lemma}
  \label{lem:separate}
  Let $F$ be a finite lower set and suppose $\varphi\in\HP$.  If
  \begin{equation*}
    M_\varphi=\begin{pmatrix} ([1]-\varphi \varphi^*)(x,y)
    \end{pmatrix}_{x,y\in F}
    \notin \mathcal C_F,
  \end{equation*}
  then there exists a cyclic unital representation $\mu:P(F)\to
  B(\mathcal H)$ such that $\|\mu(\psi)\|\le 1$ for all $\psi\in\PrF$,
  but $\|\mu(\varphi)\|>1$.
\end{lemma}

\begin{proof}
  By Theorem \ref{thm:closedcone} the cone $\mathcal C_F$ is closed
  (in the set of $F\times F$ matrices $M(F)$).  As a consequence of
  the Hahn-Banach Theorem (see, for example, \S12.F of
  \cite{MR0410335}), there is a linear functional $\lambda$ on $M(F)$
  such that $\lambda$ is nonnegative on $\mathcal C_F$ and
  $\lambda(M_\varphi) < 0$.  As $\|M_\phi\|+M_\phi \in \mathcal C_F$
  by Lemma~\ref{lem:Kpos}, we have $\lambda(1) > 0$, where $1$ is the
  identity in $M(F)$.  So in particular, $\lambda$ is not identically
  zero on $\mathcal C_F$.

  Next define a scalar product on $P(F)$ by 
  \begin{equation}
    \label{eq:defineip}
    \ip{f}{g} =\lambda(fg^*).
  \end{equation}
  For ease of notation, we will simply write ``$f$'' for the
  restriction $f|_F$ of $f$ to the lower set $F$.  We then view
  $f,g\in\mathbb C^F$ as vectors so that $fg^* \in M(F)$ is the matrix
  with entries $fg^*(x,y)=f(x)g(y)^*$.  Since, by Lemma
  \ref{lem:Kpos}, the cone $\mathcal C_F$ contains all positive
  matrices and $\lambda$ is non-negative on $\mathcal C_F$, the form
  in equation (\ref{eq:defineip}) is positive semi-definite.

  Mod out by the kernel and let $q(f)$ denote the image of $f$ in the
  quotient.  (Since the space is finite dimensional there is no need
  to complete to get a Hilbert space.)  The resulting Hilbert space,
  which we call $\mathcal H$, is nontrivial.  In particular,
  $q(\delta^F)\neq 0$.  To see this, first note that $[1]\in \mathcal
  C_F$, so $\lambda([1]) \geq 0$.  By assumption $\lambda(([1]-\varphi
  \varphi^*)) < 0$, which implies $\lambda(\varphi \varphi^*) > 0$.
  Since finite products of the test functions restricted to $F$ span
  $P(F)$, which is finite dimensional, we can write $\varphi =
  \sum_{k=0}^n c_k\xi_k$, for some finite collection of finite
  products of test functions $\{\xi_k\}$.  Repeated use of the
  equality
  \begin{equation*}
    [1]- \xi_k\xi_j\xi_j^*\xi_k^* = ([1]- \xi_k\xi_k^*) + \xi_k([1]-
    \xi_j\xi_j^*)\xi_k^*
  \end{equation*}
  and Lemma~\ref{lem:conjcone} shows that $[1]- \xi\xi^*$ is in
  $\mathcal C_F$, and so $\lambda([1]- \xi\xi^*) \geq 0$, for any
  finite product of test functions $\xi$.  By the Cauchy-Schwarz
  inequality, for any $j,k$, $|\lambda(\xi_j \xi_k^*)| \leq
  \lambda(\xi_j\xi_j^*) \lambda(\xi_k\xi_k^*)$, and so if for all $k$,
  $\lambda(\xi_k\xi_k^*)$ were zero, we would have $\lambda(\varphi
  \varphi^*) = 0$.  Hence there is some product of test functions
  $\xi$ such that $\lambda(\xi\xi^*) > 0$.  Consequently $\lambda([1])
  > 0$, and so $\|q(\delta^F)\| > 0$.

  Let $\mu$ be the right regular representation of $P(F)$ on $\mathcal
  H$.  That is, $\mu(g)q(f)=q(f\s g)$ --- provided of
  course that it is well defined.  If $\psi\in\Psi$, then because of
  the definition of $\mathcal C_F$,
  \begin{equation*}
    \|q(f)\|^2-  \| q(f\s\psi)\|^2=\lambda (f^* \s ([1]-\psi\psi^*)\sh
    f) \ge 0,
  \end{equation*}
  where the inequality follows from Lemma \ref{lem:conjcone}.  Thus,
  $\mu(\psi)$ is well defined and since $\PrF$ generates $P(F)$, $\mu$
  is well defined.

  Clearly $\mu$ is cyclic with cyclic vector $q(\delta^F)$.  Finally,
  \begin{equation*}
    \|q(\delta^F)\|^2 - \|\mu(\varphi)q(\delta^F)\|^2 =
    \lambda([1]-\varphi\varphi^*) < 0
  \end{equation*}
  so that $\| \mu(\varphi)\|>1$.
\end{proof}

\subsection{A compact set}
\label{subsec:compactset}

Fix $\varphi:G\to \mathbb C$ and a collection of test functions $\TT$.
For $F\subset G$ a finite lower set, let
\begin{equation*}
  \Phi_F = \{\Gamma\in M(F,\CT^*)^+ : ([1]-\varphi\varphi^*)(x,y) =
  (\Gamma\sh ([1]-EE^*))(x,y) \text{ for } x,y\in F\}.
\end{equation*}
The set $\Phi_F$ is naturally identified with a subset of the
product of $\CT^*$ with itself $|F|^2$ times.

\begin{lemma}
\label{lem:compact-set}
  The set $\Phi_F$ is compact.
\end{lemma}

\begin{proof}
  Let $\Gamma_\alpha$ be a net in $\Phi_F$.  Arguing as in the proof
  of Theorem \ref{thm:closedcone}, we find each $\Gamma_\alpha(x,x)$
  is a bounded net and thus each $\Gamma_{\alpha}(x,y)$ is also a
  bounded net.  By weak-$*$ compactness of the unit ball in $\CT^*$
  there exists a $\Gamma$ and subnet $\Gamma_\beta$ of $\Gamma_\alpha$
  so that for each $x,y\in F$, the net $\Gamma_{\beta}(x,y)$ converges
  to $\Gamma(x,y)$.
\end{proof}


\section{Proof of the realization theorem, Theorem \ref{thm:main}}
\label{sec:proof}

\subsection{Proof of (i) implies (iiF)}

Suppose that (iiF) does not hold.  In this case there exists a finite
lower set $F\subset G$ so that the matrix
\begin{equation*}
  M_\varphi = \begin{pmatrix} ([1]-\varphi\varphi^*)(x,y)
  \end{pmatrix}_{x,y\in F}
\end{equation*}
is not in the cone
\begin{equation*}
  \mathcal C_F=\{ \begin{pmatrix} \Gamma\sh ([1]-EE^*)
  \end{pmatrix}_{x,y\in F}: \Gamma \in M(F,\CT^*)^+\}.
\end{equation*}
Lemma \ref{lem:separate} produces a representation $\mu:P(F)\to
B(\mathcal H)$ so that $\|\mu(\psi)\|\le 1$ for all $\psi \in \PrF$,
but $\|\mu(\pi_F(\varphi))\|>1$.  Lemma \ref{lem:reptok} now implies
$\|\varphi\|>1$.

\subsection{Proof of (iiF) implies (iiG)}

The proof here uses Kurosh's Theorem and in much the same way as in
\cite{MR1882259}.

The hypothesis is that for every finite lower set $F\subset G$,
$\Phi_F$, as defined in Subsection \ref{subsec:compactset} is not
empty.  The result in that section is that $\Phi_F$ is compact.  For a
finite lower set $F$ contained in a lower set $H$, define $ \pi_F^H
:\Phi_H \to \Phi_F$ by
\begin{equation*}
  \pi_F^H (\Gamma)  =\Gamma|_{F\times F}.
\end{equation*}
Thus, with $\mathcal F$ equal to the collection of all finite lower
subsets of $G$ partially ordered by inclusion, the triple
$(\Phi_G,\pi_F^G,\mathcal F)$ is an inverse limit of nonempty compact
spaces.  Consequently, by Kurosh's Theorem (\cite{MR1882259}, p.~30),
for each $F\in\mathcal F$ there is a $\Gamma_F \in \Phi_F$ so that
whenever $F,H\in\mathcal F$ and $F\subset H$,
\begin{equation}
  \label{eq:consistent}
  \pi_F^H (\Gamma_H)  =\Gamma_F.
\end{equation}

Define $\Gamma:G\times G\to \CT^*$ by $\Gamma(x,y)=\Gamma_{F}(x,y)$
where $F\in\mathcal F$ is any lower set so that $x,y\in F$.  This is
well defined by the relation in equation (\ref{eq:consistent}).  If
$F$ is any finite lower set and $f:F\to \CT$ is any function, then
\begin{equation*}
  \sum_{x,y\in F} \Gamma(x,y)(f(x)f(y)^*) =\sum_{x,y\in
    F}\Gamma_F(x,y)(f(x)f(y)^*) \geq 0
\end{equation*}
since $\Gamma_F\in M(F,\CT^*)^+$.  Any finite subset of $G$ is
contained in a finite lower set, and so it follows that $\Gamma$ is
positive.

\subsection{Proof of (iiG) implies (iii)}

Let $\Gamma$ denote the positive kernel in (iiG).  Apply Proposition
\ref{prop:factorization} to find $\mathcal E$, $L:G\to B(\CT, \mathcal
E)$, and $\rho:\CT\to B(\mathcal E)$ as in the conclusion of the
proposition.
  
Rewrite condition (iiG) as
\begin{alignat}{1}
  \label{eq:iso1}
  \nonumber [1](x,y)+ (\Gamma\sh(EE^*))(x,y)
  =\, & (\varphi\varphi^*)(x,y) + (\Gamma\sh[1])(x,y)\\
  \nonumber \Updownarrow & \\
  \nonumber [1](x,y)+
  {\displaystyle{\sum_{pq=x}\sum_{rs=y}}}\Gamma(q,s)(E(p)E(r)^*)
  =\, & \varphi (x)
  \varphi (y)^* + \displaystyle{\sum_{pq=x}\sum_{rs=y}}
  \Gamma(q,s)(\delta(p)\delta(r)^*)\\
  \nonumber \Updownarrow & \\
  \nonumber \delta(x)\delta(y)^* + 
  {\displaystyle{\sum_{pq=x}\sum_{rs=y}}}\ip{L(q)E(p)1}{L(s)E(r)1}
  =\, & \varphi (x) \varphi (y)^* + \displaystyle{\sum_{pq=x}\sum_{rs=y}}
  \ip{L(q)\delta(p)1}{L(s)\delta(r)1}\\
  \nonumber \Updownarrow & \\
  \delta(x)\delta(y)^* + \ip{(\rho (E) \s L)(x)1}{(\rho (E) \s L)(y)1}
  =\, & \varphi (x) \varphi (y)^* + \ip{(\rho (\delta) \s L)(x)1}{(\rho
    (\delta) \s L)(y)1},
\end{alignat}
where $1$ is the identity in $\CT$.  We have used the intertwining
relation between $L$ and $\rho$ from
Proposition~\ref{prop:factorization}.  Notice that in doing so the
$\sth$-product is replaced by the $\st$-product.

From here the remainder of the proof is the standard lurking isometry
argument.

Let $\mathcal E_d$ denote finite linear combinations of
\begin{equation*}
  \begin{pmatrix} (\rho (E) \s L)(x)1 \\ \delta(x) \end{pmatrix} \in
  \begin{matrix} \mathcal E \\ \oplus \\ \mathbb C \end{matrix}
\end{equation*}
and let $\mathcal E_r$ denote finite linear combinations of
\begin{equation*}
  \begin{pmatrix} (\rho(\delta)\s L)(x)1 \\ \varphi(x) \end{pmatrix}
  \in \begin{matrix} \mathcal E \\ \oplus \\ \mathbb C \end{matrix}.
\end{equation*}

Define $V:\mathcal E_d\to \mathcal E_r$ by
\begin{equation*}
  V\begin{pmatrix} (\rho (E) \s L)(x)1 \\ \delta(x) \end{pmatrix} =
  \begin{pmatrix} (\rho(\delta)\s L)(x)1 \\ \varphi(x) \end{pmatrix}
\end{equation*}
and extending by linearity.  Equation (\ref{eq:iso1}) implies
\begin{equation*}
  \begin{split}
    {\left\| \sum c_j \begin{pmatrix} (\rho (E) \s L)(x_j)1 \\
          \delta(x_j) \end{pmatrix} \right\|}^2
    &\,= \sum_{j,\ell} c_j c_\ell^* ([1]_{x_j,x_\ell}+\ip{(\rho (E) \s
      L)(x_j)1)}{(\rho (E) \s L)(x_\ell)1)} ) \\
    &\,= \sum_{j,\ell} c_j c_\ell^* (\varphi(x_j) \varphi(x_\ell)^* +
    \ip{(\rho(\delta)\s L)(x_j)1}{\rho(\delta)\s L)(x_\ell)1} )\\
    &\,={\left\| \sum c_j \begin{pmatrix} \rho(\delta)\s L)(x_j)1 \\
          \varphi(x_j) \end{pmatrix} \right\|}^2
  \end{split}
\end{equation*}
which shows simultaneously that $V$ is well defined and an isometry.
Thus $V$ (the lurking isometry) extends to an isometry from the
closure of $\mathcal E_d$ to the closure of $\mathcal E_r$.  There
exists a Hilbert space $\mathcal H$ containing $\mathcal E$ and a
unitary map
\begin{equation*}
  U: \begin{matrix} \mathcal H \\ \oplus \\\mathbb C \end{matrix} 
  \to  \begin{matrix} \mathcal H \\ \oplus \\ \mathbb C \end{matrix}
\end{equation*}
so that $U$ restricted to $\mathcal E_d$ is $V$; i.e.,
$U\gamma=V\gamma$ for $\gamma\in\mathcal E_d.$

Write
\begin{equation}
  \label{eq:iso3}
  U=\begin{pmatrix} A & B\\ C & D \end{pmatrix}
\end{equation}
with respect to the decomposition $\mathcal H\oplus \mathbb C$.  In
particular,
\begin{equation*}
  \begin{pmatrix} A & B\\ C & D \end{pmatrix} \begin{pmatrix}
    (\rho (E) \s L)(x)1 \\ \delta(x) \end{pmatrix}= \begin{pmatrix}
    (\rho(\delta)\s L)(x)1 \\ \varphi(x) \end{pmatrix}
\end{equation*}
which gives the system of equations
\begin{eqnarray}
  \label{eq:sys1}
  A (\rho (E) \s L)(x)1 +B\delta(x) &=& (\rho(\delta)\s L)(x)1 \\
  \nonumber C (\rho (E) \s L)(x)1+D\delta(x) &=& \varphi(x).
\end{eqnarray}
From the first equation in (\ref{eq:sys1}) we have
\begin{equation}
  \label{eq:sys2}
  L(x)1 = ((\rho(\delta)-A\rho(E))^{-1_\st}\s (B\delta))(x),
\end{equation}
where the inverse is with respect to the $\st$-product.  Plugging
this into the second equation of \eqref{eq:sys1} gives
\begin{equation}
  \varphi(x)=D\delta(x) +C(\rho (E) \s (\rho(\delta) -A\rho(E))^{-1_\st}\s
  (B\delta))(x),
\end{equation}
which, using the fact that $\rho$ is unital, can be written
\begin{equation}
  \label{eq:transferrep}
  \varphi(x)=D\delta(x) +C(\rho (E) \s (\delta -A\rho(E))^{-1_\st}\s
  (B\delta))(x),
\end{equation}
as desired.

\subsection{Proof of (iii) implies (i)}

Suppose $\varphi=W_\Sigma$ as in equation (\ref{eq:transfer})
(equivalently, equation (\ref{eq:transferrep}) above).  We want to
show that $\mathbf{k} - \varphi \s \mathbf{k} \sh \varphi^* \geq 0$
for all $\mathbf{k}\in \mathcal K_\Psi$.  First, factor
$\mathbf{k}(x,y) = k_x k_y^*$.  To make the notation consistent with
that used in the calculations below, we write $k(x)$ for $k_x$.  We
compute the (four) terms in $(\mathbf{k}-\varphi \s \mathbf{k} \sh
\varphi^*)(x,y)$, using the identities implied by $U$ being unitary
and the equality $(\delta\s k )(k^*\sh \delta) = \delta \s \mathbf{k}
\sh \delta = \mathbf{k}$.  Recall that for functions $f$ and $g$,
$(f\s g)^* = g^*\sh f^*$.

To begin with, $CC^* = 1-DD^*$ and so we have $D \mathbf{k}(x,y) D^* =
\mathbf{k}(x,y)- C\mathbf{k}(x,y)C^*$.  Hence
\begin{equation*}
  \begin{split}
    &\mathbf{k}(x,y) - D\mathbf{k}(x,y) D^*= (C(\delta \s k))(x)((k^*\sh
    \delta)C^*)(y)\\ 
    =\,& \left(C(\delta-\rho(E)A)^{-1_\st}\s
      (\delta-\rho(E)A)\s k\right)(x) \left(k^*\sh
      (\delta-\rho(E)A)^* \sh (\delta-\rho(E)A)^{*-1_\st}C^*\right)(y) .
  \end{split}
\end{equation*}

For the next few terms it is useful to observe that
\begin{equation*}
  \rho(E)\s (\delta - A\rho(E))^{-1_\st} = (\delta -
  \rho(E)A)^{-1_\st} \s \rho(E),
\end{equation*}
or equivalently,
\begin{equation*}
  (\delta - \rho(E)A)^{-1_\st} \s \rho(E) \s (\delta -
  A\rho(E))= \rho(E),
\end{equation*}
which follows from
\begin{equation*}
  \rho(E) \s (\delta - A\rho(E))= \rho(E)\s \delta -
  \rho(E)\s A\rho(E) = (\delta - \rho(E)A) \s \rho(E).
\end{equation*}

The second term we consider is
\begin{equation*}
  \begin{split}
    & \left(C\rho(E)\s (\delta-A\rho(E))^{-1_\st}\s
      (B\delta)\s k\right)(x) (k^*\sh \delta^*D^*)(y) \\
    =\,& \left(C (\delta-\rho(E)A)^{-1_\st}\s\rho(E)\s
      (Bk)\right)(x) (k^*D^*)(y) \\
    =\,& \sum_{pqr=x} C(\delta-\rho(E)A)^{-1_\st}(p)\rho(E(q))
    Bk(r)k(y)^*D^* \\
    =\,& \sum_{pqr=x} C(\delta-\rho(E)A)^{-1_\st}(p)\rho(E(q))
    \mathbf{k}(r,y) BD^* \\
    =\,& -\sum_{pqr=x} C(\delta-\rho(E)A)^{-1_\st}(p)\rho(E(q))
    A \mathbf{k}(r,y) C^* \\
    =\,& -\sum_{pqr=x} C(\delta-\rho(E)A)^{-1_\st}(p)\rho(E(q))
    A k(r)k(y)^* C^* \\
    =\,& - \left(C(\delta-\rho(E)A)^{-1_\st}\s \rho(E)
    \s Ak\right)(x) \left(k^* \sh
      (\delta-\rho(E)A)^* \sh (\delta-\rho(E)A)^{*-1_\st}C^*\right)(y).
    \\ 
  \end{split}
\end{equation*}

For the third term,
\begin{equation*}
  \begin{split}
    & (D\delta\s k)(x) \left(C\rho(E)\s (\delta-A\rho(E))^{-1_\st}\s
      (B\delta)\s k\right)(y)^* \\
    =\,& (D k)(x) \left(C(\delta-\rho(E)A)^{-1_\st}\s\rho(E)\s
      (B\delta)\s k\right)(y)^* \\
    =\,& (D k)(x) \left((k^* B^*)\sh \rho(E)^* \sh
      (\delta-\rho(E)A)^{*-1_\st} C^*\right)(y) \\
    =\,& -(C k)(x) \left( A^* k^* \sh \rho(E)^* \sh
      (\delta-\rho(E)A)^{*-1_\st} C^*\right)(y) \\
    =\,& -\left(C(\delta-\rho(E)A)^{-1_\st}\s
      (\delta-\rho(E)A)\s k\right)(x) \left( A^*k^*\sh \rho(E)^* \sh 
    (\delta-\rho(E)A)^{*-1_\st}C^*\right)(y).  \\
  \end{split}
\end{equation*}

Finally, the last term is
\begin{equation*}
  \begin{split}
    & \left(C\rho(E)\s (\delta-A\rho(E))^{-1_\st}\s
      (B\delta)\s k\right)(x)  \left(C\rho(E)\s
      (\delta-A\rho(E))^{-1_\st}\s (B\delta)\s k\right)(y)^* \\
    =\,& \left(C (\delta-\rho(E) A)^{-1_\st} \s \rho(E) \s
      (B k)\right)(x)  \left(C (\delta-\rho(E)A)^{-1_\st}
      \s \rho(E) \s (B k)\right)(y)^* \\
    =\,& \left(C (\delta-\rho(E) A)^{-1_\st} \s \rho(E) \s
      (B k)\right)(x)  \left((k^*B^*)\sh \rho(E)^* \sh
      (\delta-\rho(E)A)^{*-1_\st}C^*\right)(y) \\
    =\,& \left(C (\delta-\rho(E)A)^{-1_\st}\s \rho(E)\s k\right)(x) (1-AA^*) 
      \left(k^*\sh\rho(E)^*\sh (\delta-\rho(E)A)^{*-1_\st}
        C^*\right)(y). \\ 
  \end{split}
\end{equation*}

Putting them together we have
\begin{equation*}
  \begin{split}
    &\mathbf{k}-\varphi \s \mathbf{k} \sh \varphi^* \\
    =\,& \left(C(\delta-\rho(E)A)^{-1_\st}\s
      (\delta-\rho(E)A) \s k\right) \left(k^* \sh
      (\delta-\rho(E)A)^* \sh (\delta-\rho(E)A)^{*-1_\st}C^*\right) \\
    &\quad
    - \left(C(\delta-\rho(E)A)^{-1_\st}\s\rho(E)
    \s Ak\right) \left(k^* \sh
      (\delta-\rho(E)A)^* \sh (\delta-\rho(E)A)^{*-1_\st}C^*\right) \\
    &\quad
    - \left(C(\delta-\rho(E)A)^{-1_\st}\s
      (\delta-\rho(E)A) \s k\right) \left( A^*k^* \sh \rho(E)^* \sh
      (\delta-\rho(E)A)^{*-1_\st}C^*\right) \\
    &\quad
    + \left(C (\delta-\rho(E)A)^{-1_\st}\s \rho(E) \s k\right) (1-AA^*) 
      \left(k^* \sh \rho(E)^*\sh (\delta-\rho(E)A)^{*-1_\st}
        C^*\right)(y), \\ 
  \end{split}
\end{equation*}
which after a bit of algebra is seen to simplify to
\begin{equation*}
  \mathbf{k}-\varphi \s \mathbf{k} \sh \varphi^* =
  C(\delta-\rho(E)A)^{-1_\st}\s
  (\mathbf{k}-\rho(E) \s \mathbf{k} \sh \rho(E)^*) \sh
  (\delta-\rho(E)A)^{*-1_\st}C^*.
\end{equation*}
Given a finite lower set $F\subset G$, the matrix
\begin{equation*}
  P=\begin{pmatrix} (\mathbf{k}-E \s \mathbf{k} \sh E^*)(x,y)
  \end{pmatrix}_{x,y\in F} \in M(F,\CT)
\end{equation*}
is a positive since its value at $\psi\in \Psi$ is
\begin{equation*}
  P(\psi) = \begin{pmatrix} (\mathbf{k} - \psi \s \mathbf{k} \sh
    \psi^*)(x,y) \end{pmatrix}.
\end{equation*}
Consequently $\mathbf{k}-\varphi \s \mathbf{k} \sh \varphi^* \geq 0$
over $F$, which completes the proof.


\section{Agler-Jury-Pick interpolation}

We now turn to the proof of Theorem \ref{thm:main-interpolate}.
Condition (i) implies condition (ii) simply by the definition of the
norm on $\HP$.

If condition (iii) does not hold, then an argument just as in the
proof of Theorem \ref{thm:main} produces a kernel
$\mathbf{k}\in\mathcal K_\Psi$ so that the relevant kernel on $F$ is
not positive.  Hence (ii) implies (iii).

To prove that (iii) implies (i), first argue along the lines of the
proof of (iiG) implies (iii) in Theorem \ref{thm:main}, but work with
the finite set $F$ in place of $G$.  Next verify that the transfer
function $W_\Sigma$ so constructed and defined on all of $G$ satisfies
$W_\Sigma(x)=f(x)$ for $x\in F$ (since we worked with $F$).  The
implication (iii) implies (i) in Theorem \ref{thm:main} now says that
$\|W_\Sigma\|\le 1$.\qed

\medskip

This leads to the following, which is reminiscent of results on
left tangential Nevanlinna-Pick interpolation.

\begin{theorem}
  \label{thm:NP-interpolation}
  Let $F$ be a finite lower set in a semigroupoid $G$.  Suppose
  $w(a),z(a) \in \mathbb C$, $a\in F$ are given.  Then there is a
  function $\varphi\in \HP$ with $\|\varphi\|_{\HP} \le 1$ such that
  \begin{equation*}
    (\varphi\s z)(a) = w(a), \qquad \text{for all }a\in F,
  \end{equation*}
  if and only if
  \begin{equation*}
    (z^*z - w^*w) \s \mathbf{k} \geq 0, \qquad \text{for all
    }\mathbf{k}\in K_\Psi.
  \end{equation*}
\end{theorem}

\begin{proof}
  If $w =\varphi \s z$ with $\varphi\in \HP$ and $\|\varphi\|_{\HP} \le
  1$, then for all $\mathbf{k}\in K_\Psi$,
  \begin{equation*}
    \begin{split}
      (z^*z - w^*w) \s \mathbf{k} &= \left(z^* \sh ([1] -\varphi^* \varphi) \s
        z\right) \s \mathbf{k} \\
      &= ([1] -\varphi^* \varphi) \s (z^*z) \s \mathbf{k} \\
      &= ([1] -\varphi^* \varphi) \s (z^* \s \mathbf{k} \sh z) \\
      & \geq 0,
    \end{split}
  \end{equation*}
  since by Lemma~\ref{lem:conjcone}, $z^* \s \mathbf{k} \sh z \in K_\Psi$.

  Now suppose $(zz^* - ww^*) \s \mathbf{k} \geq 0$ for all
  $\mathbf{k}\in K_\Psi$.  Begin by assuming that $z$ is an invertible
  function with respect to the $\st$~product (i.e., $z(a)$ is
  invertible for all $a\in F_e$).  Set $f = w\s z^{-1_\st}$ on $F$.
  Then restricting to $F$,
  \begin{equation*}
     0\leq (z^*z - w^*w) \s \mathbf{k} = \left(z^* \sh ([1] -f^* f) \s
        z\right) \s \mathbf{k}  = ([1] -f^* f) \s (z^* \s \mathbf{k} \sh z)
  \end{equation*}
  for all $\mathbf{k}\in K_\Psi$.  Again by Lemma~\ref{lem:conjcone},
  $z^{-1_\st\,*} \s \mathbf{k} \sh z^{-1_\st}$ is in $K_\Psi$ if
  $\mathbf{k}\in K_\Psi$.  Hence $([1] -f^* f) \s \mathbf{k} \geq 0$
  on $F\times F$ for all $\mathbf{k}\in K_\Psi$, and so by
  Theorem~\ref{thm:main-interpolate} $f$ extends to $\varphi\in \HP$
  with $\|\varphi\|_{\HP} \le 1$ such that $([1] -\varphi^* \varphi)
  \s \mathbf{k} \geq 0$ on $G\times G$.  Pad $z$ with zeros to make it
  a function in $H^2(\mathbf{k})$ for all $\mathbf{k}\in K_\psi$, and
  set $w = \varphi \s z$ (which agrees with the original definition of
  $w$ on the lower set $F$).

  If $z$ is not $\st$-invertible, then $z(a) = 0$ for some $a \in F_e$.
  This means that $\{z(a): a \in F\}$ (where $z(a)$ is identified with
  the vector with this value in the $a^{\mathrm{th}}$ position and
  zero elsewhere) is not a basis for $\mathbb C^F$.  Choose a vector
  $g$ with $g(a) = 1$ for each $a$ in $F_e$ where $z(a) = 0$ and zero
  otherwise.  Fix $\epsilon > 0$.  Let $z'$ be $g$ normalized so that
  $\mathrm{Re}\,\ip{z}{z'} \geq 0$ and $\|z'\| < \epsilon$.  Then for
  $z_\epsilon = z+z'$,
  \begin{equation*}
    (z_\epsilon^* z_\epsilon -w^*w)\s \mathbf{k} \geq (z^*z + {z'}^*z'
    - w^*w)\s \mathbf{k} \geq (z^*z - w^*w)\s \mathbf{k} \geq 0,
  \end{equation*}
  and $z_\epsilon$ is invertible, so we obtain by the last paragraph a
  corresponding $f_\epsilon$ for which $([1] -f_\epsilon^* f_\epsilon)
  \s \mathbf{k} \geq 0$ for all $\mathbf{k}\in K_\Psi$ and $f_\epsilon
  \s z_\epsilon = w$.  Since we are on a finite dimensional space, the
  sequence $\{f_{1/n}\}_{n=1,2,\ldots}$ converges to some $f\in P(F)$
  with $([1] -f^* f) \s \mathbf{k} \geq 0$ for all $\mathbf{k}\in
  K_\Psi$.  Also $z_{1/n} \longrightarrow z$.  Consequently, $f\s z =
  w$.
\end{proof}

A right tangential problem could very easily be formulated and solved.
One way to do this would be to replace ``$\st$'' with ``$\hat\s$'' at
appropriate points in the left interpolation theorem and proof, and
then take adjoints.  The details are left to the interested reader.

Finally note that taking $z = \delta_F$ and $w=f$ in the last
theorem recovers the first two equivalences in
Theorem~\ref{thm:main-interpolate}.


\section{Examples}
\label{sec:examples}

\subsection{The classical examples}

View $\mathbb D$ as a Pick semigroupoid.  The partial multiplication
is trivial and so each $z\in\mathbb D$ is idempotent.  Take
$\Psi=\{z\}$ ($z$ meaning here the identity function) as the
collection of test functions.  The Agler-Jury-Pick interpolation
theorem in this case is Pick interpolation.

Choose $G=\mathbb N$ with the usual semigroup(oid) structure.  Let
$\Psi=\{z\}$, where by $z$ we mean the function $z:\mathbb N\to
\mathbb C$ given by $z(j)=0$ if $j\ne 1$ and $z(1)=1$ (we think of
$z(j)$ as the derivatives of $z$ at $0$).  In this case
Agler-Jury-Pick interpolation is Carath\'eodory-Fej\'er interpolation.

For mixed Agler-Pick and Carath\'eodory-Fej\'er choose $G=\mathbb D
\times \mathbb N$ with the semigroupoid structure,
\begin{equation*}
  (z,n)(w,m)=\begin{cases} \text{ is not  defined } & \text{ if } z\ne
    w \\  (z,n+m) & \mbox{ if } z=w \end{cases}
\end{equation*}
and let $\Psi=\{z\}$ denote the function $z(w,0)=w$, $z(w,1)=1$ and
$z(w,m)=0$ for $m\ge 2$.

\subsection{Agler-Pick interpolation on an annulus}

Let $\mathbb A$ denote an annulus $\{q<|z|<1\}$, viewed as a Pick
semigroupoid.

There is a family of analytic functions $\psi:\mathbb A\to \mathbb D$
which are unimodular on the boundary of $\mathbb A$ and have precisely
two zeros in $\mathbb A$ (counting with multiplicity), normalized by
$\psi(\sqrt{q})=0$ and $\psi(1)=1$.  If $\varphi$ is any other
analytic function on $\mathbb A$ which is unimodular on the boundary
and has exactly two zeros (counting with multiplicity), then there is
a M\"obius map $m$ from the disk onto the disk such that $m\circ
\varphi \in\Psi$.  There is a canonical parameterization of $\Psi$ by
the unit circle.

\begin{theorem}
  The collection $\Psi$ is a family of test functions for $\mathbb A$
  and the norm in $\HP$ is the same as the norm on $H^\infty (\mathbb
  A)$.  Moreover, no proper subset of $\Psi$ is a set of test
  functions which gives the norm of $H^\infty(\mathbb A)$.

  In the case of Agler-Pick interpolation $($on a finite set $F\subset
  \mathbb A)$, the realization formula for a solution is in terms of a
  single positive measure on the unit circle.
\end{theorem}

Look for the details of this example in the forthcoming paper
\cite{DMnext}.

\subsection{Carath\'eodory interpolation kernels}
\label{subsec:cara-int-kernels}

Let $\mathbb N$ denote the natural numbers with the usual
semigroup(oid) structure.  A kernel $\mathbf{k}$ on $\mathbb N$ is a
Carath\'eodory interpolation kernel \cite{MR1134687} provided (by way
of normalization) $\mathbf{k}(0,0)=1$, $\mathbf{k}(0,n)=0$ for $n>0$,
and
\begin{equation*}
  b=[1]-\mathbf{k}^{-1_\st}
\end{equation*}
is positive.

For illustrative purposes, suppose $b$ has finite rank $d$ and so
factors as $b=B^*B$, where $B:\mathbb N \to (\mathbb C^d)^*$.
Although $B$ is not scalar-valued, $[1]-B(a)B(b)^*=[1]-b$ is scalar
and moreover,
\begin{equation*}
  ([1]-BB^*)\s \mathbf{k}=([1]-b)\s \mathbf{k}= k^{-1_\st}\s k=[1]\ge 0.
\end{equation*}
Choosing $\Psi=\{B\}$, it turns out that $\varphi \in \HP$ and
$\|\varphi \|_{\HP}\le 1$ if and only if $([1]-\varphi \varphi^*)\s
\mathbf{k}$ is positive.

\subsection{NP kernels and Arveson-Arias-Popescu space}

The situation for Nevanlinna-Pick (NP) kernels is similar to that for
Carath\'eodory kernels. In particular, it requires a version of our
results for vector valued test functions.

As a particular example, consider the semigroup $\mathbb N^g$ with the
(single) vector valued test function $Z=\begin{pmatrix} z_1 & z_2 &
  \cdots & z_g\end{pmatrix}^T$.  This pair $(\mathbb N^g,Z)$ gives
rise to symmetric Fock space; i.e., the space of multipliers of the
space of analytic functions on the unit ball in $\mathbb C^g$ with
reproducing kernel $\mathbf{k}(z,w)=(1-\ip{z}{w} )^{-1}$ studied by
Arveson (\cite{MR1668582,MR2054985}, in the commutative case) and by
Arias and Popescu (\cite{MR1749679,MR1675377}, in both the commutative
case and the noncommutative case discussed in the next subsection).

\subsection{Noncommutative Toeplitz algebras}
\label{subsec:nonc-toepl-algebr}

The following have been considered in the context of Nevanlinna-Pick
and Carath\'eodory-Fej\'er interpolation by Davidson and Pitts
\cite{MR1627901} and Arias and Popescu \cite{MR1749679}, as well as
by Popescu in \citep{MR1652931,MR1348353}.

Let $\mathfrak F=\mathfrak F_g$ denote the free monoid on the $g$
letters $\{x_1,\dots,x_g\}$.  Let $\psi_j:\mathfrak F\to \mathbb C$
denote the function $\psi_j(x_j)=1$ and $\psi(w)=0$ if $w$ is any word
other than $x_j$.  The matrix $\sT(\psi_j)$ is a (truncated) shift on
Fock space.

Given a word $w=x_{j_1}x_{j_2}\cdots x_{j_n}$, let
\begin{equation*}
  \psi^{w_\st}=\psi_{j_1}\s \psi_{j_2}\s \cdots \s \psi_{j_n}.
\end{equation*}
Since $\psi^{w_\st}(v)=1$ if $w=v$ and $0$ otherwise, it follows that
if $F$ is any finite subset of $\mathfrak F$, then $P(F)$ contains all
functions on $F$.

Let $\psi=\begin{pmatrix} \psi_1 & \cdots & \psi_g \end{pmatrix}^T$ and
consider $\psi$ as a (single) test function.  We calculate
\begin{equation*}
     s(x,y)=  ([1]-\psi^* \psi)^{-1_\st}
\end{equation*}
where $s$ is the Toeplitz kernel ($s(x,y)=1$ of $x=y$ and $s(x,y)=0$
if $x\ne y$). Then if $\mathbf{k}$ is any kernel for which
$([1]-\psi^* \psi)\s \mathbf{k}=Q$ is positive, we have $s\s Q =
\mathbf{k}$. It follows that $\|\varphi\|\le 1$ if and only if the
kernel
\begin{equation*}
 \mathfrak F\times \mathfrak F \ni (x,y)\mapsto  
    ([1]-\varphi \varphi^*)\s s(x,y) 
\end{equation*}
is positive, if and only if $\|\sT(\varphi)\| \le 1$.  The versions of
the Nevanlinna-Pick theorem considered in the papers cited above
coincide with Theorem~\ref{thm:NP-interpolation}, while
Carath\'eodory-Fej\'er interpolation is given by
Theorem~\ref{thm:main-interpolate}.

As a final remark, note that each $\sT(\psi_j)$ is an isometry and
\begin{equation*}
  \sum \sT(\varphi)\sT(\varphi)^*= P_{\emptyset} \geq 0.
\end{equation*}
Here $P_{\emptyset}$ is the projection onto the span of the vacuum
vector $\emptyset$ in the Fock space.

\subsection{The Polydisk}

The semigroupoid $\mathbb N^g$ (the $g$-fold product of the
nonnegative integers) with the set of test functions $z_j$, the
characteristic function of $e_j$ the vector with $1$ in the $j$-th
entry and $0$ elsewhere, gives rise to the Schur-Agler class of the
polydisk $\mathbb D^g$ returning us to the introduction
and~\cite{MR1207393}.

\subsection{Semigroupoid algebras of Power and Kribs}

Kribs and Power \citep{MR2076898,MR2063759} consider a
generalization of the noncommutative Toeplitz algebras which they term
a free semigroupoid algebra.  Order arises from the assumption of
freeness, the resulting semigroupoid is cancellative, and there is a
representation (related to our Toeplitz representation on
characteristic functions $\chi_a$) in terms of partial isometries and
projections.

A notion of a generalized Fock space is developed, which is simply the
Hilbert space with orthonormal basis labelled by the elements of the
semigroupoid.  The algebras of interest in these papers are obtained
from the weak operator topology closure of the algebras generated by
the left regular representations (i.e., the projections and partial
isometries mentioned above).

The algebras are closely related to those in the present paper when
$G$ is a semigroupoid in this more restrictive sense and the
collection of test functions consists of the characteristic functions
of non-idempotent elements from the first stratum (to use our
terminology).

It is assumed that for every idempotent $e\in G$, there is a
non-idempotent $a$ such that $ae$ is defined.  Let $G_1$ be the first
(left) stratum in $G$, and assume that this set is countable.  Then
$G$ is generated by $G_1$, in the sense that if $x$ is in the
$n^{\mathrm{th}}$ stratum, then $x=ay$, where $y$ is in the
$(n-1)^{\mathrm{st}}$ stratum and $a$ is in the first stratum.  Let
$P$ have the property that
\begin{equation*}
  P(x,y) =
  \begin{cases}
    1 & x=y \text{ and }x,y\notin G_e,\\
    0 & \text{otherwise,}
  \end{cases}
\end{equation*}
and $\tilde s = [1]+P$.  Clearly $\tilde s$ is invertible.  
Now mimic the proof in Section~\ref{subsec:nonc-toepl-algebr} by
letting $\psi_j(x_j)=1$ if $x_j\in G_1$.  It is not difficult to
verify that for $\psi=\begin{pmatrix} \psi_1 & \dots & \psi_g
\end{pmatrix}^T$,
\begin{equation*}
     \tilde s(x,y)=  ([1]-\psi^* \psi)^{-1_\st},
\end{equation*}
and so just as in that subsection, if $\mathbf{k}$ is any kernel for
which $([1]-\psi^* \psi)\s \mathbf{k}=Q$ is positive, $\tilde s\s Q =
\mathbf{k}$. It follows that the statements $\|\varphi\|\le 1$,
\begin{equation*}
 \mathfrak F\times \mathfrak F \ni (x,y)\mapsto  
    ([1]-\varphi \varphi^*)\s \tilde s(x,y) 
\end{equation*}
positive, and $\|\sT(\varphi)\| \le 1$ are all equivalent.

A number of interesting algebras can be generated in this manner,
including the noncommutative Toeplitz algebras above and the norm
closed semicrossed product ${\mathbb C}^n \times_\beta^\sigma {\mathbb
  Z}_+$ \cite{MR2063759}.  Indeed, the condition of being freely
generated can be replaced by our more general conditions for a
semigroupoid (again assuming though that for every idempotent $e\in
G$, there is a non-idempotent $a$ such that $ae$ is defined).  Our
results allow for interpolation in all of these algebras.

\bibliographystyle{plain}
\bibliography{semigroupoids}

\end{document}